\documentclass{tran-l}

\usepackage{amssymb}

\usepackage{graphicx}

\usepackage{amsfonts}
\usepackage{amsmath}
\usepackage{amscd}
\usepackage{algorithmic}
\usepackage{comment}
\usepackage{color}
\usepackage{array,amsbsy}
\usepackage{enumerate}
\usepackage[bottom]{footmisc}
\usepackage[all,cmtip]{xy}

\def\publname{}
\copyrightinfo{}{}
\def\ISSN{}
\def\volinfo{}
\def\pageinfo{}
\def\PII{}

\newtheorem{theorem}{Theorem}[section]
\newtheorem{corollary}[theorem]{Corollary}
\newtheorem{lemma}[theorem]{Lemma}

\theoremstyle{definition}
\newtheorem{definition}[theorem]{Definition}

\theoremstyle{remark}

\numberwithin{equation}{section}

\newenvironment{packed_enum}{
    \begin{enumerate}
        \setlength{\itemsep}{1pt}
        \setlength{\parskip}{0pt}
        \setlength{\parsep}{0pt}
}{\end{enumerate}}

\newenvironment{packed_itemize}{
    \begin{itemize}
        \setlength{\itemsep}{1pt}
        \setlength{\parskip}{0pt}
        \setlength{\parsep}{0pt}
}{\end{itemize}}

\newcommand{\heading}[1]{\smallskip\par\noindent{\bf #1}}

\def\claimqed{\hfill$\diamond$}

  \def\calC{{\mathcal C}}

\def\calM{{\mathcal M}}  
\def\calP{{\mathcal P}}  \def\calR{{\mathcal R}}

\def\cP{\hbox{\rm \sffamily P}}
\def\cNP{\hbox{\rm \sffamily NP}}
\def\cGI{\hbox{\rm \sffamily GI}}

\def\O{\mathcal{O}{}}

\def\int{\hbox{\rm \sffamily INT}}
\def\circle{\hbox{\rm \sffamily CIRCLE}}
\def\comp{\hbox{\rm \sffamily COMP}}
\def\cocomp{\hbox{\rm \sffamily co-COMP}}

\def\uint{\hbox{\rm \sffamily UNIT\ INT}}
\def\ca{\hbox{\rm \sffamily CIRCULAR-ARC}}
\def\chor{\hbox{\rm \sffamily CHOR}}
\def\fun{\hbox{\rm \sffamily FUN}}
\def\perm{\hbox{\rm \sffamily PERM}}
\def\caterpillar{\hbox{\rm \sffamily CATERPILLAR}}

\def\planar{\hbox{\rm \sffamily PLANAR}}
\def\tree{\hbox{\rm \sffamily TREE}}

\def\clawfree{\hbox{\rm \sffamily CLAW-FREE}}
\def\cobip{\hbox{\rm \sffamily co-BIP}}

\def\ifa{\hbox{\rm \sffamily IFA}}
\def\graphs#1{\hbox{\rm \sffamily $#1$-GRAPH}}
\def\dim#1{\hbox{\rm \sffamily $#1$-DIM}}
\def\trapezoid{\hbox{\rm \sffamily TRAPEZOID}}
\def\bipperm{\hbox{\rm \sffamily BIP PERM}}

\def\bolddim#1{\hbox{\rm \bf \sffamily $\boldsymbol{#1}$-DIM}}

\def\Aut{{\rm Aut}}

\def\Stab{{\rm Stab}}

\def\Autcc{{\rm Aut(\text{\rm connected}\ \circle)}}

\def\dm{{\rm dim}}

\def\im{{\rm Im}}
\def\id{{\rm id}}
\def\sec{{\rm sec}}

\def\twin{\thicksim_{TW}}
\def\Ker{{\rm Ker}}
\def\orient{\mathfrak{to}}

\def\sym{\mathbb{S}}
\def\dih{\mathbb{D}}
\def\cyc{\mathbb{Z}}

\def\Rep{\mathfrak{Rep}}
\def\Repfac{\mathfrak{Rep} / \mathord{\sim}}
\def\toord{\mathord{\to}}
\def\overtoord{\overline{\raisebox{0pt}[1.75mm]{$\toord$}}}

\newtheorem{problem}{Problem}
\newtheorem{conjecture}{Conjecture}

\begin{document}

\title{Automorphism Groups of Geometrically Represented Graphs}


\author{Pavel Klav\'{i}k}
\address{Computer Science Institute, Charles University in Prague, Czech Republic}
\email{klavik@iuuk.mff.cuni.cz, zeman@iuuk.mff.cuni.cz}
\thanks{The conference version of this paper appeared at STACS 2015~\cite{kz}.  For an interactive
graphical map of the results, see \url{http://pavel.klavik.cz/orgpad/geom_aut_groups.html}
(supported for Google Chrome).  Supported by CE-ITI (P202/12/G061 of GA\v{C}R) and Charles
University as GAUK 196213.}

\author{Peter Zeman}

\subjclass[2010]{Primary 05C62, 08A35, 20D45}

\date{}

\dedicatory{}

\begin{abstract}
We describe a technique to determine the automorphism group of a geometrically represented graph, by
understanding the structure of the induced action on all geometric representations. Using this, we
characterize automorphism groups of interval, permutation and circle graphs. 
We combine techniques from group theory (products, homomorphisms, actions) with data structures from
computer science (PQ-trees, split trees, modular trees) that encode all geometric representations.

We prove that interval graphs have the same automorphism groups as trees, and for a given interval
graph, we construct a tree with the same automorphism group which answers a question of Hanlon
[Trans.~Amer.~Math.~Soc 272(2), 1982]. For permutation and circle graphs, we give an inductive
characterization by semidirect and wreath products. We also prove that every abstract group can be
realized by the automorphism group of a comparability graph/poset of the dimension at most four.
\end{abstract}

\maketitle

\section{Introduction} \label{sec:introduction}

The study of symmetries of geometrical objects is an ancient topic in mathematics and its precise
formulation led to group theory. Symmetries play an important role in many distinct areas. In 1846,
Galois used symmetries of the roots of a polynomial in order to characterize polynomials which are
solvable by radicals.

\heading{Automorphism Groups of Graphs.}
The symmetries of a graph $X$ are described by its automorphism group $\Aut(X)$. Every automorphism
is a permutation of the vertices which preserves adjacencies and non-adjacencies.
Frucht~\cite{frucht1939herstellung} proved that every finite group is isomorphic to $\Aut(X)$ of
some graph $X$. General algebraic, combinatorial and topological structures can be encoded by
(possibly infinite) graphs~\cite{hedrlin} while preserving automorphism groups.

Most graphs are asymmetric, i.e., have only the trivial automorphism~\cite{er63}.  However, many
mathematical results rely on highly symmetrical objects.  Automorphism groups are important for
studying large regular objects, since their symmetries allow one to simplify and understand these
objects. 


\begin{definition}
For a graph class $\calC$, let $\Aut(\calC) = \bigl\{G : X \in \calC,\ G \cong \Aut(X)\bigr\}$.
The class $\calC$ is called \emph{universal} if every abstract finite group is contained in
$\Aut(\calC)$, and \emph{non-universal} otherwise.
\end{definition}

\newpage

In 1869, Jordan~\cite{jordan1869assemblages} gave a characterization for the class of trees (\tree):

\begin{theorem}[Jordan~\cite{jordan1869assemblages}]\label{thm:jordan_trees}
The class $\Aut(\tree)$ is defined inductively as follows:
\begin{packed_enum}
\item[\rm (a)] $\{1\} \in \Aut(\tree)$.
\item[\rm (b)] If $G_1, G_2 \in \Aut(\tree)$, then $G_1 \times G_2 \in \Aut(\tree)$.
\item[\rm (c)] If $G \in \Aut(\tree)$, then $G \wr \sym_n \in \Aut(\tree)$.
\end{packed_enum}
\end{theorem}

\noindent The direct product in (b) constructs the automorphisms that act independently on non-isomorphic
subtrees and the wreath product in (c) constructs the automorphisms that permute isomorphic
subtrees.

\heading{Graph Isomorphism Problem.} This famous problem asks whether two input graphs $X$ and $Y$
are the same up to a relabeling. This problem is obviously in \cNP, and not known to be
polynomially-solvable or \cNP-complete. Aside integer factorization, this is a prime candidate for
an intermediate problem with the complexity between \cP\ and \cNP-complete.  It belongs to the low
hierarchy of \cNP~\cite{schoning1988graph}, which implies that it is unlikely \cNP-complete. (Unless
the polynomial-time hierarchy collapses to its second level.) The graph isomorphism problem is known
to be polynomially solvable for the classes of graphs with bounded degree~\cite{luks1982isomorphism}
and with excluded topological subgraphs~\cite{grohe_marx}. The graph isomorphism problem is the
following fundamental mathematical question: given two mathematical structure, can we test their
isomorphism in some more constructive way than by guessing a mapping and verifying that it is an
isomorphism.  

The graph isomorphism problem is closely related to computing generators of an automorphism group.
Assuming $X$ and $Y$ are connected, we can test $X \cong Y$ by computing generators of $\Aut(X
\mathbin{\dot\cup} Y)$ and checking whether there exists a generator which swaps $X$ and $Y$.
For the converse relation, Mathon~\cite{mathon1979note} proved that generators of the automorphism
group can be computed using $\O(n^4)$ instances of graph isomorphism. Compared to graph isomorphism,
automorphism groups of restricted graph classes are much less understood.

\heading{Geometric Representations.} In this paper, we study automorphism groups of geometrically
represented graphs.  The main question we address is how the geometry influences their automorphism groups.
For instance, the geometry of a sphere translates to 3-connected planar graphs which have unique
embeddings~\cite{whitney}. Thus, their automorphism groups are so called spherical groups which are
the automorphism groups of tilings of a sphere. For general planar graphs (\planar), the automorphism
groups are more complex and they were described by Babai~\cite{babai1972automorphism} and in more
details in~\cite{kn} by semidirect products of spherical and symmetric groups.

We focus on intersection representations. An \emph{intersection representation} $\calR$ of a graph
$X$ is a collection $\{R_v : v \in V(X)\}$ such that $uv \in E(X)$ if and only if $R_u \cap R_v \ne
\emptyset$; the intersections encode the edges. To get nice graph classes, one typically restricts
the sets $R_v$ to particular classes of geometrical objects; for an overview, see the classical
books~\cite{agt,egr}. We show that a well-understood structure of all intersection representations
allows one to determine the automorphism group.

\begin{figure}[t!]
\centering
\includegraphics[scale=0.8]{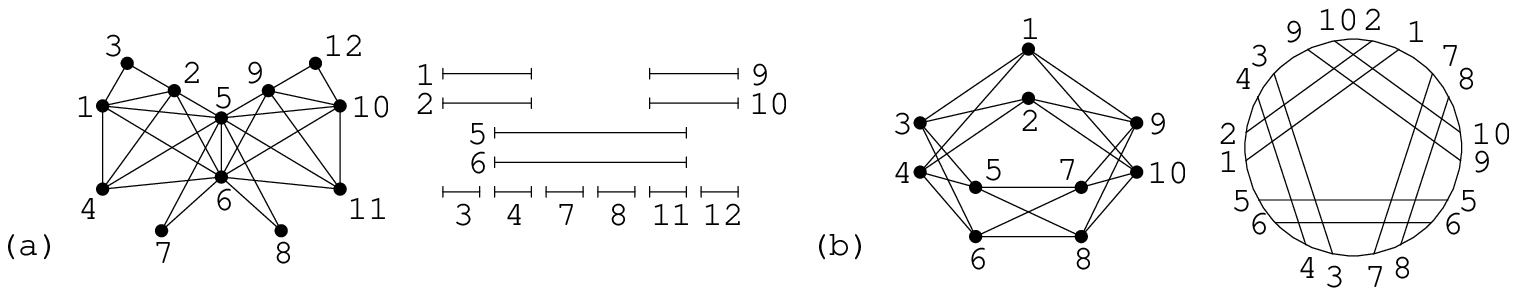}
\caption{(a) An interval graph and one of its interval representations. (b) A circle graph and one
of its circle representations.}
\label{fig:interval_graph_rep_circle_graph_rep}
\end{figure}

\heading{Interval Graphs.}
In an \emph{interval representation} of a graph, each set $R_v$ is a closed interval of the real
line.  A graph is an \emph{interval graph} if it has an interval representation; see
Fig.~\ref{fig:interval_graph_rep_circle_graph_rep}a.  A graph is a \emph{unit interval
graph} if it has an interval representation with each interval of the length one.  We denote these
classes by \int\ and \uint, respectively.  \emph{Caterpillars} (\caterpillar) are trees with
all leaves attached to a central path; we have $\caterpillar = \int \cap \tree$.

\begin{theorem} \label{thm:aut_groups_interval}
The following equalities hold:
\begin{packed_enum}
\item[\rm (i)] $\Aut(\int) = \Aut(\tree)$,
\item[\rm (ii)] $\Aut(\text{\rm connected}\ \uint) = \Aut(\caterpillar)$,
\end{packed_enum}
\end{theorem}

Concerning (i), this equality is not well known. It was stated by Hanlon~\cite{hanlon1982counting}
without a proof in the conclusion of his paper from 1982 on enumeration of interval graphs. Our
structural analysis is based on PQ-trees~\cite{PQ_trees} which describe all interval representations
of an interval graph. It explains this equality and further solves an open problem of Hanlon: for a
given interval graph, to construct a tree with the same automorphism group. Without PQ-trees, this
equality is surprising since these classes are very different. Caterpillars which form their
intersection have very restricted automorphism groups (see
Lemma~\ref{prop:automorphisms_caterpillars}). The result (ii) follows from the known
properties of unit interval graphs and our understanding of $\Aut(\int)$. 

\heading{Circle Graphs.}
In a \emph{circle representation}, each $R_v$ is a chord of a circle. A graph is a circle graph
(\circle) if it has a circle representation; see
Fig.~\ref{fig:interval_graph_rep_circle_graph_rep}b.

\begin{theorem} \label{thm:aut_groups_circle}
Let $\Sigma$ be the class of groups defined inductively as follows:
\begin{packed_itemize}
\item[\rm (a)] $\{1\} \in \Sigma$.
\item[\rm (b)] If $G_1, G_2 \in \Sigma$, then $G_1
        \times G_2 \in \Sigma$.
\item[\rm (c)] If $G \in \Sigma$, then $G \wr
        \sym_n \in \Sigma$.
\item[\rm (d)] If $G_1, G_2, G_3, G_4 \in \Sigma$,
        then $(G_1^4 \times G_2^2 \times G_3^2 \times G_4^2) \rtimes \cyc_2^2 \in
        \Sigma$.
\end{packed_itemize}
Then $\Aut(\text{\rm connected}\ \circle)$ consists of the following groups:
\begin{packed_itemize}
\item If $G \in \Sigma$, then $G \wr \cyc_n \in \Aut(\text{\rm connected}\ \circle)$.
\item If $G_1, G_2 \in \Sigma$, then $(G_1^n \times G_2^{2n}) \rtimes \dih_n \in \Aut(\text{\rm
connected}\ \circle)$.
\end{packed_itemize}
\end{theorem}

The automorphism group of a disconnected circle graph can be easily determined using
Theorem~\ref{thm:aut_disconnected}.  We are not aware of any previous results on the automorphism
groups of circle graphs. We use split trees describing all representations of circle graphs.
The class $\Sigma$ consists of the stabilizers of vertices in connected circle graphs and
$\Aut(\tree) \subsetneq \Sigma$.  

\begin{figure}[t!]
\centering
\includegraphics[scale=0.8]{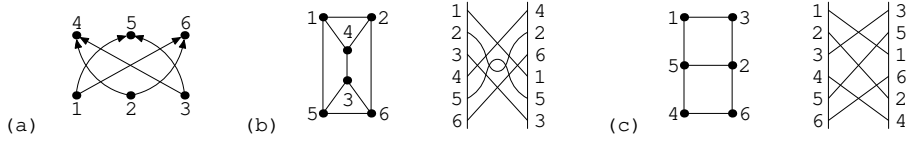}
\caption{(a) A comparability graph with a transitive orientation. (b)
A function graph and one of its representations. (c) A permutation graph and
one of its representations.}
\label{fig:comp_fun_perm_example}
\end{figure}

\heading{Comparability Graphs.} A \emph{comparability graph} is derived from a poset by removing the
orientation of the edges.  Alternatively, every comparability graph $X$ can be transitively
oriented: if $x \to y$ and $y \to z$, then $xz \in E(X)$ and $x \to z$; see
Fig~\ref{fig:comp_fun_perm_example}a.  This class was first studied by
Gallai~\cite{gallai1967transitiv} and we denote it by \comp.

An important structural parameter of a poset $P$ is its \emph{Dushnik-Miller
dimension}~\cite{dushnik_miller}.  It is the least number of linear orderings $L_1,\dots,L_k$ such
that $P = L_1 \cap \cdots \cap L_k$.  (For a finite poset $P$, its dimension is always finite since
$P$ is the intersection of all its linear extensions.) Similarly, we define the \emph{dimension} of
a comparability graph $X$, denoted by $\dm(X)$, as the dimension of any transitive orientation of
$X$. (Every transitive orientation has the same dimension; see Section~\ref{sec:aut_comp}.) By
$\dim{k}$, we denote the subclass consisting of all comparability graphs $X$ with $\dm(X) \le k$. We
get the following infinite hierarchy of graph classes:
$$\dim1 \subsetneq \dim2 \subsetneq \dim3 \subsetneq \dim 4 \subsetneq \cdots \subsetneq \comp.$$
For instance, \cite{planar_incidence_3dim} proves that the bipartite graph of the incidence between
the vertices and the edges of a planar graph always belongs to \dim3.

Surprisingly, comparability graphs are related to intersection graphs, namely to function and
permutation graphs.  \emph{Function graphs} (\fun) are intersection graphs of continuous real-valued
function on the interval $[0,1]$.  \emph{Permutation graphs} (\perm) are function graphs which can
be represented by linear functions called \emph{segments}~\cite{permutation_graphs}; see
Fig.~\ref{fig:comp_fun_perm_example}b and~c. We have $\fun =
\cocomp$~\cite{golumbic1983comparability} and $\perm = \comp \cap \cocomp =
\dim2$~\cite{even1972permutation}, where $\cocomp$ are the complements of comparability graphs.

Since \dim1 consists of all complete graphs, $\Aut(\dim1) = \{\sym_n : n \in \mathbb N\}$. The
automorphism groups of $\dim2 = \perm$ are the following:

\begin{theorem} \label{thm:aut_groups_perm}
The class $\Aut(\perm)$ is described inductively as follows:
\begin{packed_itemize}
\item[\rm (a)] $\{1\} \in \Aut(\perm)$,
\item[\rm (b)] If $G_1,G_2 \in \Aut(\perm)$, then $G_1 \times G_2 \in \Aut(\perm)$.
\item[\rm (c)] If $G \in \Aut(\perm)$, then $G \wr \sym_n \in \Aut(\perm)$.
\item[\rm (d)] If $G_1,G_2,G_3 \in \Aut(\perm)$, then $(G_1^4 \times G_2^2 \times G_3^2) \rtimes
\cyc_2^2 \in \Aut(\perm)$.
\end{packed_itemize}
\end{theorem}

In comparison to Theorem~\ref{thm:jordan_trees}, there is the additional operation (d) which 
shows that $\Aut(\tree) \subsetneq \Aut(\perm)$. Geometrically, the group $\cyc_2^2$ in (d)
corresponds to the horizontal and vertical reflections of a symmetric permutation representation.
Notice that it is more restrictive than the operation (d) in Theorem~\ref{thm:aut_groups_circle}.
Our result also easily gives the automorphism groups of \emph{bipartite permutation graphs}
(\bipperm), in particular $\Aut(\caterpillar) \subsetneq \Aut(\bipperm) \subsetneq \Aut(\perm)$.

\begin{corollary} \label{cor:aut_groups_bipperm}
The class $\Aut({\rm connected}\ \bipperm)$ consists of all abstract groups $G_1$, $G_1 \wr \cyc_2
\times G_2 \times G_3$, and $(G_1^4 \times G_2^2) \rtimes \cyc_2^2$, where $G_1$ is a direct
product of symmetric groups, and $G_2$ and $G_3$ are symmetric groups.
\end{corollary}

Comparability graphs are universal since they contain bipartite graphs; we can orient all edges from
one part to the other. Since the automorphism group is preserved by complementation, $\fun =
\cocomp$ implies that also function graphs are universal. In Section~\ref{sec:comparability_graphs},
we explain the universality of \fun\ and \comp\ in more detail using the induced action on the set
of all transitive orientations. Similarly posets are known to be
universal~\cite{universal_posets,universal_posets2}.

It is well-known that bipartite graphs have arbitrarily large dimensions: the \emph{crown graph},
which is $K_{n,n}$ without a matching, has the dimension $n$.  We give a construction which encodes any
graph $X$ into a comparability graph $Y$ with $\dm(Y) \le 4$, while preserving the automorphism
group.

\begin{theorem}\label{thm:kdim_aut_groups_and_gi}
For every $k \ge 4$, the class $\dim{k}$ is universal and its graph isomorphism is $\cGI$-complete.
The same holds for posets of the dimension $k$.
\end{theorem}

Yannakakis~\cite{yannakakis1982complexity} proved that recognizing \dim3\ is $\cNP$-complete by a
reduction from $3$-coloring. For a graph $X$, a comparability graph $Y$ is constructed with several vertices
representing each element of $V(X) \cup E(X)$. It is proved that $\dm(Y) = 3$ if and
only if $X$ is $3$-colorable. Unfortunately, the automorphisms of $X$ are lost in $Y$ since it
depends on the labels of $V(X)$ and $E(X)$ and $Y$ contains some additional edges according to these
labels. We describe a simple and completely different construction which achieves only the dimension
4, but preserves the automorphism group: for a given graph $X$, we create $Y$ by replacing each edge
with a path of length eight. However, it is non-trivial to show that $Y \in \dim4$, and the
constructed four linear orderings are inspired by~\cite{yannakakis1982complexity}. A different
construction follows from~\cite{grid_dim4,grid_GI_complete}.

\heading{Related Graph Classes.}
Theorems~\ref{thm:aut_groups_interval}, \ref{thm:aut_groups_circle} and~\ref{thm:aut_groups_perm}
and Corollary~\ref{cor:aut_groups_bipperm} state that \int, \uint, \circle, \perm, and \bipperm\ are
non-universal.  Figure~\ref{fig:graph_classes} shows that their superclasses are already universal.

\emph{Trapezoidal graphs} (\trapezoid) are intersection graphs of trapezoids between two parallel
lines and they have universal automorphism groups~\cite{trapezoidal_gi_complete}.  \emph{Claw-free
graphs} (\clawfree) are graphs with no induced $K_{1,3}$. Roberts~\cite{roberts1969indifference}
proved that $\uint = \clawfree \cap \int$. The \emph{complements of bipartite graphs} (\cobip) are
claw-free and universal.  \emph{Chordal graphs} (\chor) are intersection graphs of subtrees of
trees. They contain no induced cycles of length four or more and naturally generalize interval
graphs. Chordal graphs are universal~\cite{lueker1979linear}.  \emph{Interval filament graphs}
(\ifa) are intersection graphs of the following sets.  For every $R_u$, we choose an interval
$[a,b]$ and $R_u$ is a continuous function $[a,b] \to \mathbb{R}$ such that $R_u(a) = R_u(b) = 0$
and $R_u(x) > 0$ for $x \in (a, b)$.

\begin{figure}[b!]
\centering
\includegraphics[scale=0.75]{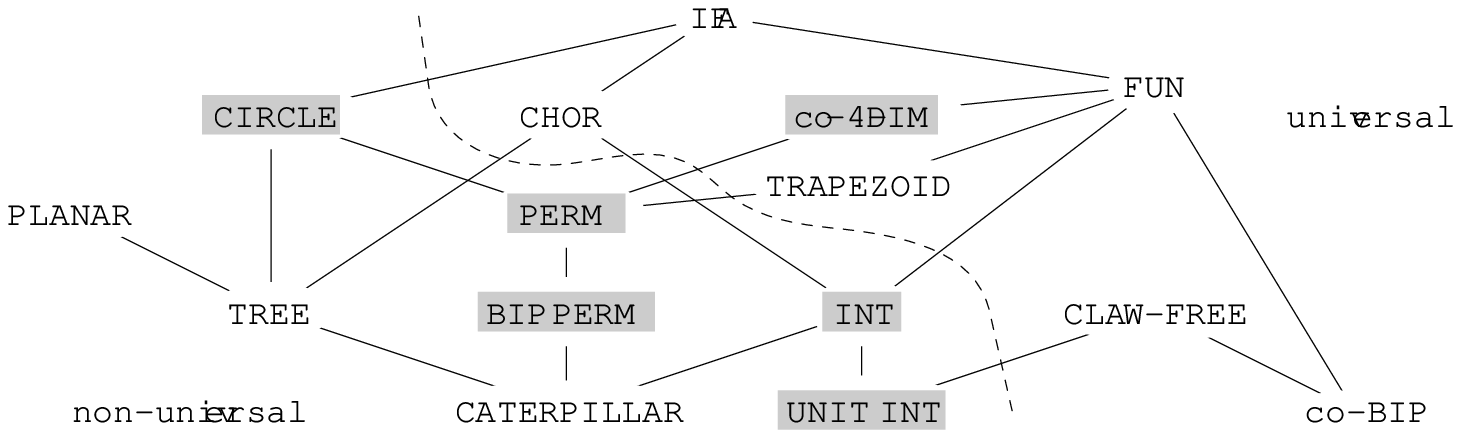}
\caption{The inclusions between the considered graph classes. We characterize the automorphism groups of
the classes in gray.}
\label{fig:graph_classes}
\end{figure}

\heading{Outline.}
In Section~\ref{sec:preliminaries}, we introduce notation and group products.  In
Section~\ref{sec:action}, we explain our general technique for determining the automorphism group
from the geometric structure of all representations, and relate it to map theory. We describe the
automorphism groups of interval and unit interval graphs in Section~\ref{sec:interval_graphs}, of
circle graphs in Section~\ref{sec:circle_graphs}, and of permutation and bipartite permutation
graphs in Section~\ref{sec:comparability_graphs}.  Our results are constructive and lead to
polynomial-time algorithms computing automorphism groups of these graph classes; see
Section~\ref{sec:algorithms}.  We conclude with several open problems.

\section{Preliminaries} \label{sec:preliminaries}

We use $X$ and $Y$ for graphs, $M$, $T$ and $S$ for trees and $G$ and $H$ for groups. The
vertices and edges of $X$ are $V(X)$ and $E(X)$. For $A \subseteq V(X)$, we denote by $X[A]$ the
subgraph induced by $A$, and for $x \in V(X)$, the closed neighborhood of $x$ by $N[x]$.
The complement of $X$ is denoted by $\overline{X}$, clearly $\Aut(X) = \Aut(\overline{X})$.

A permutation $\pi$ of $V(X)$ is an \emph{automorphism} if $uv \in E(X) \iff \pi(u)\pi(v) \in E(X)$.
The automorphism group $\Aut(X)$ consists of all automorphisms of $X$.  We use the notation
$\sym_n$, $\dih_n$ and $\cyc_n$ for the \emph{symmetric}, \emph{dihedral} and \emph{cyclic groups}.
Note that $\dih_1 \cong \cyc_2$ and $\dih_2 \cong \cyc_2^2$ (which appears in
Theorems~\ref{thm:aut_groups_circle} and~\ref{thm:aut_groups_perm} in (d)).
An action is called \emph{semiregular} if all stabilizers are trivial.

\heading{Group Products.}
Group products allow decomposing of large groups into smaller ones.
%
Given two groups $N$ and $H$, and a group homomorphism $\varphi \colon H \to \Aut(N)$,
we can construct a new group $N \rtimes_{\varphi} H$ as the Cartesian product $N \times H$ with
the operation defined as $(n_1, h_1) \cdot (n_2,h_2) = (n_1 \cdot \varphi(h_1)(n_2), h_1\cdot
h_2)$.  The group $N\rtimes_{\varphi}H$ is called the \emph{external semidirect product of $N$ and $H$
with respect to the homomorphism $\varphi$}, and sometimes we omit the homomorphism $\varphi$ and
write $N \rtimes H$. Alternatively, $G$ is the \emph{internal semidirect product of $N$ and $H$} if
$N \unlhd G$, $H \le G$, $N \cap H$ is trivial and $\left<N \cup H\right> = G$.

Suppose that $H$ acts on $\{1,\dots,n\}$.  The \emph{wreath product} $G \wr H$ is a shorthand for
the semidirect product $G^n \rtimes_{\psi} H$ where $\psi$ is defined naturally by $\psi(\pi) =
(g_1,\dots,g_n) \mapsto (g_{\pi(1)},\dots,g_{\pi(n)})$. In the paper, we have $H$ equal $\sym_n$ or
$\cyc_n$ for which we use the natural actions on $\{1,\dots,n\}$. For more details,
see~\cite{carter2009visual,rotman1995introduction}. All semidirect products used in this paper are
\emph{generalized wreath products} of $G_1,\dots,G_k$ with $H$, in which each orbit
of the action of $H$ has assigned one group $G_i$.

\subsection{Automorphism Groups of Disconnected Graphs.}

In 1869, Jordan described the automorphism groups of disconnected graphs, in terms of the automorphism groups
of their connected components. Since a similar argument is used in several places in this paper, we
describe his proof in details. Figure~\ref{fig:two_edges_cayley} shows the automorphism group for
a graph consisting of two isomorphic components.

\begin{figure}[!b]
\centering
\includegraphics[scale=0.8]{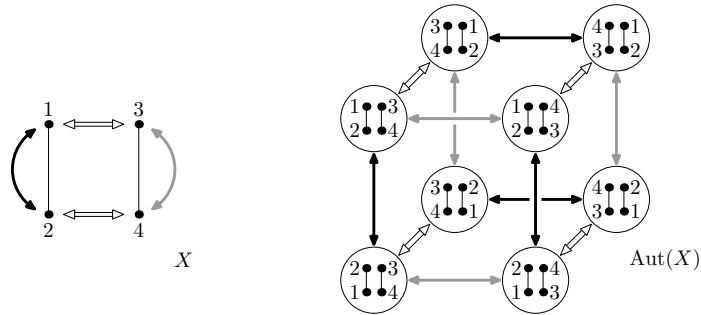}
\caption{The structure of $\Aut(X)$, generated by three involutions acting on $X$ on the left:
$\Aut(X) \cong \cyc_2^2 \rtimes \cyc_2 = \cyc_2 \wr \cyc_2$.}
\label{fig:two_edges_cayley}
\end{figure}

\begin{theorem}[Jordan~\cite{jordan1869assemblages}] \label{thm:aut_disconnected}
If $X_1, \dots, X_n$ are pairwise non-isomorphic connected graphs and $X$ is the
disjoint union of $k_i$ copies of $X_i$, then
$$\Aut(X) \cong \Aut(X_1) \wr \sym_{k_1} \times \cdots \times \Aut(X_n) \wr \sym_{k_n}.$$
\end{theorem}

\begin{proof}
Since the action of $\Aut(X)$ is independent on non-isomorphic components, it is clearly the direct
product of factors, each corresponding to the automorphism group of one isomorphism class of
components. It remains to show that if $X$ consists of $k$ isomorphic components of a connected
graph $Y$, then $\Aut(X) \cong \Aut(Y) \wr \sym_k$.

We isomorphically label the vertices of each component. Then each automorphism
$\pi \in \Aut(X)$ is a composition $\sigma \cdot \tau$ of two automorphisms: $\sigma$ maps each
component to itself, and $\tau$ permutes the components as in $\pi$ while preserving the
labeling. Therefore, the automorphisms $\sigma$ can be bijectively identified with the elements of
$\Aut(Y)^k$ and the automorphisms $\tau$ with the elements of $\sym_k$.

Let $\pi,\pi' \in \Aut(X)$.  Consider the composition $\sigma \cdot \tau \cdot \sigma' \cdot \tau'$,
we want to swap $\tau$ with $\sigma'$ and rewrite this as a composition $\sigma \cdot \hat\sigma
\cdot \hat\tau \cdot \tau$. Clearly the components are permuted in $\pi \cdot \pi'$ exactly as in
$\tau \cdot \tau'$, so $\hat\tau = \tau$. On the other hand, $\hat\sigma$ is not necessarily equal
$\sigma'$. Let $\sigma'$ be identified with the vector $(\sigma_1',\dots,\sigma_k') \in \Aut(Y)^k$.
Since $\sigma'$ is applied after $\tau$, it acts on the components permuted according to $\tau$.
Therefore $\hat\sigma$ is constructed from $\sigma'$ by permuting the coordinates of its vector by $\tau$:
$$\hat\sigma = (\sigma'_{\tau(1)},\dots,\sigma'_{\tau(k)}).$$
This is precisely the definition of the wreath product, so $\Aut(X) \cong \Aut(Y) \wr \sym_k$.
\end{proof}

\subsection{Automorphism Groups of Trees.}

Using the above, we can explain why $\Aut(\tree)$ is closed under (b) and (c):

\begin{proof}[Proof of Theorem~\ref{thm:jordan_trees} (a sketch)]
We assume that trees are rooted since the automorphism groups preserve centers.  Every
inductively defined group can be realized by a tree as follows.  For the direct product in (b), we
choose two non-isomorphic trees $T_1$ and $T_2$ with $\Aut(T_i) \cong G_i$, and attach them to a
common root. For the wreath product in (c), we take $n$ copies of a tree $T$ with $\Aut(T) \cong G$
and attach them to a common root. On the other hand, given a rooted tree, we can delete the root and
apply Theorem~\ref{thm:aut_disconnected} to the created forest of rooted trees.
\end{proof}

\section{Automorphism Groups Acting on Intersection Representations} \label{sec:action}

In this section, we describe the general technique which allows us to
geometrically understand automorphism groups of some intersection-defined graph classes.
Suppose that one wants to understand an abstract group $G$. Sometimes, it is possible interpret $G$
using a natural action on some set which is easier to understand. The action is called
\emph{faithful} if no element of $G$ belongs to all stabilizers. The structure of $G$ is captured
by a faithful action. We require that this action is ``faithful enough'', which means that the
stabilizers are simple and can be understood.

Our approach is inspired by map theory. A \emph{map} $\calM$ is a 2-cell embedding of a graph; i.e,
aside vertices and edges, it prescribes a rotational scheme for the edges incident with each vertex.
One can consider the action of $\Aut(X)$ on the set of all maps of $X$: for $\pi \in \Aut(X)$, we get
another map $\pi(\calM)$ in which the edges in the rotational schemes are permuted by $\pi$; see
Fig.~\ref{fig:maps_example}. The stabilizer of a map $\calM$, called the automorphism group
$\Aut(\calM)$, is the subgroup of $\Aut(X)$ which preserves/reflects the rotational schemes. Unlike
$\Aut(X)$, we know that $\Aut(\calM)$ is always small (since $\Aut(\calM)$ acts semiregularly on
flags) and can be efficiently determined.  The action of $\Aut(X)$ describes morphisms between
different maps and in general can be very complicated. Using this approach, the automorphism groups
of planar graphs can be characterized~\cite{babai1972automorphism,kn}.

\begin{figure}[t!]
\centering
\includegraphics[scale=0.8]{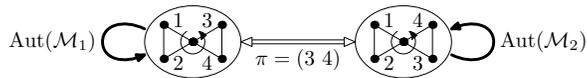}
\caption{There are two different maps, depicted with the action of $\Aut(X)$. The stabilizers
$\Aut(\calM_i) \cong \cyc_2^2$ are normal subgroups. The remaining automorphisms morph one map into
the other, for instance $\pi$ transposing $2$ and $3$. We have $\Aut(X) \cong
\cyc_2^2 \rtimes \cyc_2$.}
\label{fig:maps_example}
\end{figure}

\heading{The Induced Action.} For a graph $X$, we denote by $\Rep$ the set of all its
(interval, circle, etc.) intersection representations. An automorphism $\pi \in \Aut(X)$ creates
from $\calR \in \Rep$ another representation $\calR'$ such that $R'_{\pi(u)} = R_u$; so $\pi$ swaps
the labels of the sets of $\calR$. We denote $\calR'$ as $\pi(\calR)$, and $\Aut(X)$ acts on $\Rep$. 

The general set $\Rep$ is too large. Therefore, we define a suitable equivalence relation $\sim$ and
we work with $\Repfac$.  It is reasonable to assume that $\sim$ is a congruence with respect to the
action of $\Aut(X)$: for every $\calR \sim \calR'$ and $\pi \in \Aut(X)$, we have $\pi(\calR) \sim
\pi(\calR')$. We consider the induced action of $\Aut(X)$ on $\Repfac$.

The stabilizer of $\calR \in \Repfac$, denoted by $\Aut(\calR)$, describes automorphisms inside this
representation. For a nice class of intersection graphs, such as interval, circle or permutation
graphs, the stabilizers $\Aut(\calR)$ are very simple.  If it is a normal subgroup, then the
quotient $\Aut(X) / \Aut(\calR)$ describes all morphisms which change one representation in the
orbit of $\calR$ into another one. Our strategy is to understand these morphisms geometrically, for
which we use the structure of all representations, encoded for the considered classes by PQ-, split
and modular trees.

\section{Automorphism Groups of Interval Graphs} \label{sec:interval_graphs}

In this section, we prove Theorem \ref{thm:aut_groups_interval}. We introduce an MPQ-tree which
combinatorially describe all interval representations of a given interval graph.  We define its
automorphism group, which is a quotient of the automorphism group of the interval graph.
Using MPQ-trees, we derive a characterization of $\Aut(\int)$ which we prove to be equivalent to the
Jordan's characterization of $\Aut(\tree)$. Finally, we solve Hanlon's open
problem~\cite{hanlon1982counting} by constructing for a given interval graph a tree with the same
automorphism group, and we also show the converse construction.

\subsection{PQ- and MPQ-trees}

We denote the set of all maximal cliques of $X$ by $\calC(X)$.  In 1965, Fulkerson and Gross proved
the following fundamental characterization of interval graphs by orderings of maximal cliques: 

\begin{lemma}[Fulkerson and Gross~\cite{maximal_cliques}] \label{lem:fulkerson_gross}
A graph $X$ is an interval graph if and only if there exists a linear ordering $\preceq$ of $\calC(X)$
such that for every $x \in V(X)$ the maximal cliques containing $x$ appear consecutively in this
ordering.
\end{lemma}

\begin{proof}[Sketch of proof]
Let $\calR = \bigl\{R_x \colon x \in V(X)\bigr\}$ be an interval representation of $X$ and let
$\calC(X) = \{C_1, \dots, C_k\}$. By Helly's Theorem, the intersection $\bigcap_{x \in C_i} R_x$ is
non-empty, and therefore it contains a point $c_i$.  The ordering of $c_1, \dots, c_k$ from left to
right gives the ordering $\preceq$.

For the other implication, given an ordering $C_1 \preceq \cdots \preceq C_k$ of the maximal
cliques, we place points $c_1, \dots, c_k$ in this ordering on the real line. To each vertex $x$, we
assign the minimal interval $R_x$ such that $c_i \in R_x$ for all $x \in C_i$.  We obtain a valid
interval representation $\bigl\{R_x \colon x \in V(X)\bigr\}$ of $X$.
\end{proof}

An ordering $\preceq$ of $\calC(X)$ from Lemma~\ref{lem:fulkerson_gross} is called a
\emph{consecutive ordering}. Consecutive orderings of $\calC(X)$ correspond to different
interval representations of $X$.

\heading{PQ-trees.}
Booth and Lueker~\cite{PQ_trees} invented a data structure called a \emph{PQ-tree} which encodes all
consecutive orderings of an interval graph. They build this structure to construct a linear-time
algorithm for recognizing interval graphs which was a long standing open problem. PQ-trees give a
lot of insight into the structure of all interval representations, and have applications to many
problems. We use them to capture the automorphism groups of interval graphs.  

A rooted tree $T$ is a PQ-tree representing an interval graph $X$ if the following holds.
It has two types of inner nodes: \emph{P-nodes} and \emph{Q-nodes}. For every inner node, its
children are ordered from left to right. Each P-node has at least two children and each Q-node at
least three. The leaves of $T$ correspond one-to-one to $\calC(X)$. The \emph{frontier} of
$T$ is the ordering $\preceq$ of the leaves from left to right.

Two PQ-trees are \emph{equivalent} if one can be obtained from the other by a sequence of two
\emph{equivalence transformations}: (i) an arbitrary permutation of the order of the children of a
P-node, and (ii) the reversal of the order of the children of a Q-node. The consecutive orderings of
$\calC(X)$ are exactly the frontiers of the PQ-trees equivalent with $T$.  Booth and
Lueker~\cite{PQ_trees} proved the existence and uniqueness of PQ-trees (up to equivalence
transformations). Figure~\ref{fig:interval_graph_pq_tree_mpq_tree} shows an example.

For a PQ-tree $T$, we consider all sequences of equivalent transformations.  Two such sequences are
\emph{congruent} if they transform $T$ the same.  Each sequence consists of several transformations of
inner nodes, and it is easy to see that these transformation are independent. If a sequence
transforms one inner node several times, it can be replaced by a single transformation of this node.
Let $\Sigma(T)$ be the quotient of all sequences of equivalent transformations of $T$ by this
congruence. We can represent each class by a sequence which transforms each node at most once.

Observe that $\Sigma(T)$ forms a group with the concatenation as the group operation. This
group is isomorphic to a direct product of symmetric groups.  The order of $\Sigma(T)$ is equal to
the number of equivalent PQ-trees of $T$.  Let $T' = \sigma(T)$ for some $\sigma \in \Sigma(T)$.
Then $\Sigma(T') \cong \Sigma(T)$ since $\sigma' \in \Sigma(T')$ corresponds to $\sigma \sigma'
\sigma^{-1} \in \Sigma(T)$. 

\begin{figure}[t]
\centering
\includegraphics[scale=0.8]{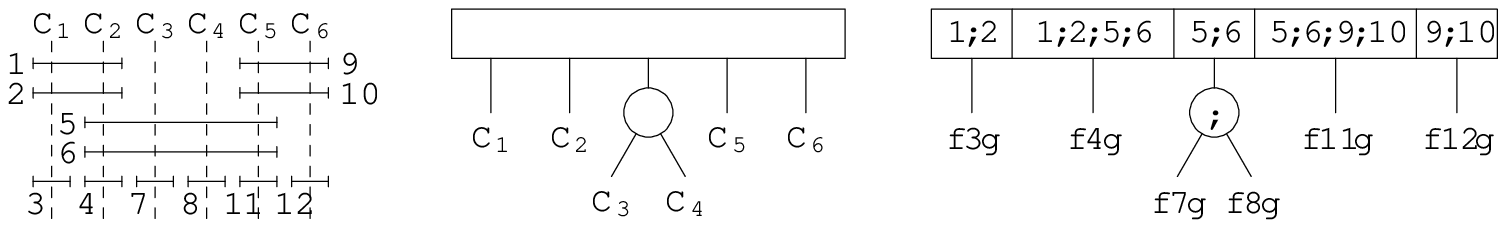}
\caption{An ordering of the maximal cliques, and the corresponding PQ-tree and MPQ-tree. The P-nodes
are denoted by circles, the Q-nodes by rectangles. There are four different consecutive orderings.}
\label{fig:interval_graph_pq_tree_mpq_tree}
\end{figure}

\heading{MPQ-trees.}
A \emph{modified PQ-tree} is created from a PQ-tree by adding information about the vertices. They
were described by Korte and M\"{o}hring~\cite{incremental_linear_int_recognition} to simplify
linear-time recognition of interval graphs. It is not widely known but the equivalent idea was used
earlier by Colbourn and Booth~\cite{autmorphism_algorithms_trees_int_planar}.

Let $T$ be a PQ-tree representing an interval graph $X$. We construct the MPQ-tree~$M$ by assigning
subsets of $V(X)$, called \emph{sections}, to the nodes of $T$; see
Fig.~\ref{fig:interval_graph_pq_tree_mpq_tree}. The leaves and the P-nodes have each assigned exactly
one section while the Q-nodes have one section per child. We assign these sections as follows:
\begin{itemize}
\item For a leaf $L$, the section $\sec(L)$ contains those vertices that are only
in the maximal clique represented by $L$, and no other maximal clique.
\item For a P-node $P$, the section $\sec(P)$ contains those vertices that are in
all maximal cliques of the subtree of $P$, and no other maximal clique.
\item For a Q-node $Q$ and its children $T_1, \dots, T_n$, the section $\sec_i(Q)$
contains those vertices that are in the maximal cliques represented by the
leaves of the subtree of $T_i$ and also some other $T_j$, but not in any other
maximal clique outside the subtree of $Q$. We put $\sec(Q) = \sec_1(Q) \cup
\cdots \cup \sec_n(Q)$.
\end{itemize}
Korte and M\"{o}hring~\cite{incremental_linear_int_recognition} proved existence of MPQ-trees and
many other properties, for instance each vertex appears in sections of exactly one node and in the
case of a Q-node in consecutive sections. Two vertices are in the same sections if and only if they
belong to precisely the same maximal cliques. Figure~\ref{fig:interval_graph_pq_tree_mpq_tree} shows
an example.

We consider the equivalence relation $\twin$ on $V(X)$ is defined as follows: $x \twin y$ if and
only if $N[x] = N[y]$. If $x \twin y$, then we say that they are \emph{twin vertices}. The
equivalence classes of $\twin$ are called \emph{twin classes}. Twin vertices can usually be ignored,
but they influence the automorphism group.  Two vertices belong to the same sections if and only if
they are twin vertices.

\subsection{Automorphisms of MPQ-trees}

For a graph $X$, the automorphism group $\Aut(X)$ induces an action on $\calC(X)$ since every
automorphism permutes the maximal cliques. If $X$ is an interval graph, then a consecutive ordering
$\preceq$ of $\calC(X)$ is permuted into another consecutive ordering $\pi(\preceq)$, so $\Aut(X)$
acts on consecutive orderings.

Suppose that an MPQ-tree $M$ representing $X$ has the frontier $\preceq$. For every automorphism
$\pi \in \Aut(X)$, there exists the unique MPQ-tree $M'$ with the frontier $\pi(\preceq)$ which is
equivalent to $M$.  We define a mapping
$$\Phi : \Aut(X) \to \Sigma(M)$$
such that $\Phi(\pi)$ is the sequence of equivalent transformations which transforms $M$ to $M'$. It
is easy to observe that $\Phi$ is a group homomorphism.

By Homomorphism Theorem, we know that $\im(\Phi) \cong \Aut(X) / \Ker(\Phi)$.  The kernel
$\Ker(\Phi)$ consists of all automorphisms which fix the maximal cliques, so they permute the
vertices inside each twin class. It follows that $\Ker(\Phi)$ is isomorphic to a direct product of
symmetric groups. So $\im(\Phi)$ almost describes $\Aut(X)$.

Two MPQ-trees $M$ and $M'$ are \emph{isomorphic} if the underlying PQ-trees are equal and
there exists a permutation $\pi$ of $V(X)$ which maps each section of $M$ to the corresponding section of
$M'$. In other words, $M$ and $M'$ are the same when ignoring the labels of the vertices in the
sections.  A sequence $\sigma \in \Sigma(M)$ is called an \emph{automorphism} of $M$ if $\sigma(M)
\cong M$; see Fig.~\ref{fig:mpq_tree_automorphism}. The automorphisms of $M$ are closed under
composition, so they form the automorphism group $\Aut(M) \le \Sigma(M)$.

\begin{figure}[b]
\centering
\includegraphics[scale=0.8]{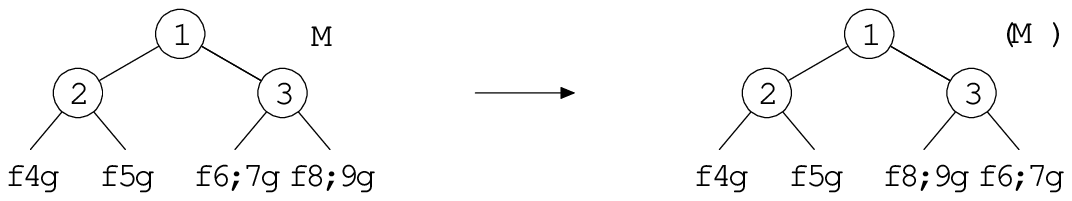}
\caption{The sequence $\sigma$, which transposes the children of the P-node with the section
$\{3\}$, is an automorphism since $\sigma(M) \cong M$. On the other, the transposition of the
children the root P-node is not an automorphism.} \label{fig:mpq_tree_automorphism}
\end{figure}

\begin{lemma} \label{lem:mpq_tree_group}
For an MPQ-tree $M$, we have $\Aut(M) = \im(\Phi)$.
\end{lemma}

\begin{proof}
Suppose that $\pi \in \Aut(X)$. The sequence $\sigma = \Phi(\pi)$ transforms $M$ into $\sigma(M)$.
It follows that $\sigma(M) \cong M$ since $\sigma(M)$ can be obtained from $M$ by permuting the
vertices in the sections by $\pi$. So $\sigma \in \Aut(M)$ and $\im(\Phi) \le \Aut(M)$.

On the other hand, suppose $\sigma \in \Aut(M)$. We know that $\sigma(M) \cong M$ and let $\pi$ be a
permutation of $V(X)$ from the definition of the isomorphism. Two vertices of $V(X)$ are adjacent if
and only if they belong to the sections of $M$ on a common path from the root. This property is
preserved in $\sigma(M)$, so $\pi \in \Aut(X)$. Each maximal clique is the union of all sections on
the path from the root to the leaf representing this clique.  Therefore the maximal cliques are
permuted by $\sigma$ the same as by $\pi$. Thus $\Phi(\pi) = \sigma$ and $\Aut(M) \le \im(\Phi)$.
\end{proof}

\begin{lemma}\label{lem:interval_graphs_semidirect}
For an MPQ-tree $M$ representing an interval graph $X$, we have $\Aut(X) \cong \Ker(\Phi) \rtimes
\Aut(M)$.
\end{lemma}

\begin{proof}
Let $\sigma \in \Aut(M)$. In the proof of Lemma~\ref{lem:mpq_tree_group}, we show that every
permutation $\pi$ from the definition of $\sigma(M) \cong M$ is an automorphism of $X$ mapped by
$\Phi$ to $\sigma$. Now, we want to choose these permutations consistently for all elements of
$\Aut(M)$ as follows. Suppose that $\id = \sigma_1,\sigma_2,\dots,\sigma_n$ be the elements of
$\Aut(M)$. We want to find $\id = \pi_1, \pi_2, \dots, \pi_n$ such that $\Phi(\pi_i) = \sigma_i$ and
if $\sigma_i\sigma_j = \sigma_k$, then $\pi_i\pi_j = \pi_k$. In other words, $H =
\{\pi_1,\dots,\pi_n\}$ is a subgroup of $\Aut(X)$ and $\Phi \restriction_H$ is an isomorphism
between $H$ and $\Aut(M) = \im(\Phi)$.

Suppose that $\pi,\pi' \in \Aut(X)$ such that $\Phi(\pi) = \Phi(\pi')$. Then $\pi$ and $\pi'$
permute the maximal cliques the same and they can only act differently on twin vertices, i.e., $\pi
\pi'^{-1} \in \Ker(\Phi)$. Suppose that $C$ is a twin class, then $\pi(C) = \pi'(C)$ but they can
map the vertices of $C$ differently. To define $\pi_1,\dots,\pi_n$, we need to define them on the
vertices of the twin classes consistently. To do so, we arbitrarily order the vertices in each twin
class. For each $\pi_i$, we know how it permutes the twin classes, suppose a twin class $C$ is
mapped to a twin class $\pi_i(C)$. Then we define $\pi_i$ on the vertices of $C$ in such a way that
the orderings are preserved.

The above construction of $H$ is correct. Since $H$ is the complementary subgroup of $\Ker(\Phi)$,
we get $\Aut(X)$ as the internal semidirect product $\Ker(\Phi) \rtimes H \cong \Ker(\Phi) \rtimes
\Aut(M)$. Our approach is similar to the proof of Theorem~\ref{thm:aut_disconnected}, and the
external semidirect product can be constructed in the same way.
\end{proof}

\subsection{The Inductive Characterization}

Let $X$ be an interval graph, represented by an MPQ-tree $M$. By
Lemma~\ref{lem:interval_graphs_semidirect}, $\Aut(X)$ can be described from $\Aut(M)$ and
$\Ker(\Phi)$. We build $\Aut(X)$ inductively using $M$, similarly as in
Theorem~\ref{thm:jordan_trees}:

\begin{proof}[Proof of Theorem~\ref{thm:aut_groups_interval}(i)]
We show that $\Aut(\int)$ is closed under (b), (c) and (d); see
Fig.~\ref{fig:interval_aut_groups_constructions}. For (b), we attach interval graphs $X_1$ and $X_2$
such that $\Aut(X_i) = G_i$ to an asymmetric interval graph. For (c), let $G \in \Aut(\int)$ and let
$Y$ be a connected interval graph with $\Aut(Y) \cong G$. We construct $X$ as the disjoint union of
$n$ copies of $Y$. For (d), we construct $X$ by attaching $X_1$ and $X_2$ to a path, where
$\Aut(X_i) = G_i$.

\begin{figure}[t!]
\centering
\includegraphics[scale=0.8]{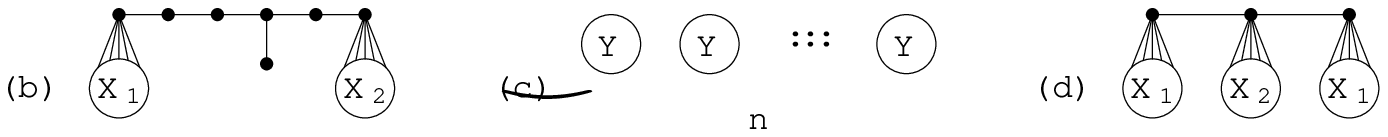}
\caption{The constructions in the proof of Theorem~\ref{thm:aut_groups_interval}(i).}
\label{fig:interval_aut_groups_constructions}
\end{figure}

For the converse, let $M$ be an MPQ-tree representing an interval graph $X$.
Let $M_1,\dots,M_k$ be the subtrees of the root of $M$ and let $X_i$ be the interval graphs induced
by the vertices of the sections of $M_i$. We want to build $\Aut(X)$ from
$\Aut(X_1),\dots,\Aut(X_k)$ using (b) to (d).

\emph{Case 1: The root is a P-node $P$.} Each sequence $\sigma \in \Aut(M)$ corresponds to
interior sequences in $\Aut(M_i)$ and some reordering $\sigma'$ of $M_1,\dots,M_k$. If
$\sigma'(M_i) = M_j$, then necessarily $X_i \cong X_j$. On each isomorphism class of
$X_1,\dots,X_k$, the permutations $\sigma'$ behave to $\Aut(X_i)$ like the permutations $\tau$ to
$\Aut(Y)$ in the proof of Theorem~\ref{thm:aut_disconnected}.  Therefore the point-wise stabilizer
of $\sec(P)$ in $\Aut(X)$ is constructed from $\Aut(X_1),\dots,\Aut(X_k)$ as in
Theorem~\ref{thm:aut_disconnected}. Since every automorphim preserves $\sec(P)$, then $\Aut(X)$ is
obtained by the direct product of the above group with the symmetric group of order $|\sec(P)|$.
Thus the operations (b) and (c) are sufficient.\footnote{Alternatively, we can show that each $X_i$
is connected and $X$ is the disjoint union of $X_1,\dots,X_k$ together with $|\sec(P)|$ vertices
attached to everything.  So Theorem~\ref{thm:aut_disconnected} directly applies.}

\emph{Case 2: The root is a Q-node $Q$.} We call $Q$ symmetric if it is transformed
by some sequence of $\Aut(M)$, and asymmetric otherwise. Let $M_1,\dots,M_k$ be its children from
left to right.  If $Q$ is asymmetric, then $\Aut(M)$ is the direct product
$\Aut(X_1),\dots,\Aut(X_k)$ together with the symmetric groups for all twin classes of $\sec(Q)$, so
it can be build using (b).  If $Q$ is symmetric, let $G_1$ is the direct product of the left part of
the children and twin classes, and $G_2$ of the middle part.  We get
$$\Aut(X) \cong (G_1^2 \times G_2) \rtimes \cyc_2 \cong G_1^2 \rtimes \cyc_2 \times G_2 \cong G_1
\wr \cyc_2 \times G_2,$$
where the wreath product with $\cyc_2$ adds the automorphisms reversing $Q$, corresponding to
reversing of vertically symmetric parts of a representation. Therefore $\Aut(X)$ can be generated
using (b) and (c).
\end{proof}

\begin{figure}[b!]
\centering
\includegraphics[scale=0.8]{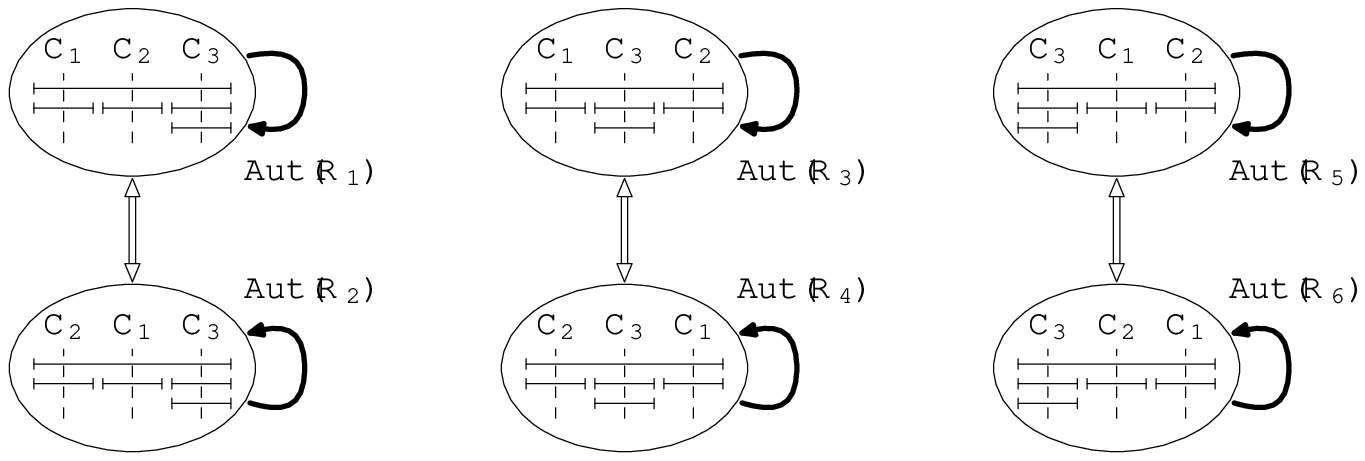}
\caption{An interval graph with six non-equivalent representation. The action of $\Aut(X)$ has three
isomorphic orbits.}
\label{fig:non-transitive_action}
\end{figure}

\subsection{The Action on Interval Representations.}

For an interval graph $X$, the set $\Rep$ consists of all assignments of closed intervals which
define $X$. It is natural to consider two interval representations equivalent if one can be
transformed into the other by continuous shifting of the endpoints of the intervals while preserving
the correctness of the representation. Then the representations of $\Repfac$ correspond to 
consecutive orderings of the maximal cliques; see Fig.~\ref{fig:non-transitive_action}
and~\ref{fig:interval_action_diagram}.

We interpret our results in terms of the action of $\Aut(X)$ on
$\Rep$. In Lemma~\ref{lem:interval_graphs_semidirect}, we proved that $\Aut(X) \cong \Ker(\Phi)
\rtimes \Aut(M)$ where $M$ is an MPQ-tree.  If an automorphism belongs to $\Aut(\calR)$, then it
fixes the ordering of the maximal cliques and it can only permute twin vertices.  Therefore
$\Aut(\calR) = \Ker(\Phi)$ since each twin class consists of equal intervals, so they can be
arbitrarily permuted without changing the representation.  Every stabilizer $\Aut(\calR)$ is the
same and every orbit of the action of $\Aut(X)$ is isomorphic, as in
Fig.~\ref{fig:non-transitive_action}.

Different orderings of the maximal cliques correspond to different reorderings of $M$. The defined
$\Aut(M) \cong \Aut(X) / \Aut(\calR)$ describes morphisms of representations belonging to one orbit
of the action of $\Aut(X)$, which are the same representations up to the labeling of the intervals;
see Fig.~\ref{fig:non-transitive_action} and Fig.~\ref{fig:interval_action_diagram}.

\begin{figure}[t!]
\centering
\includegraphics[scale=0.75]{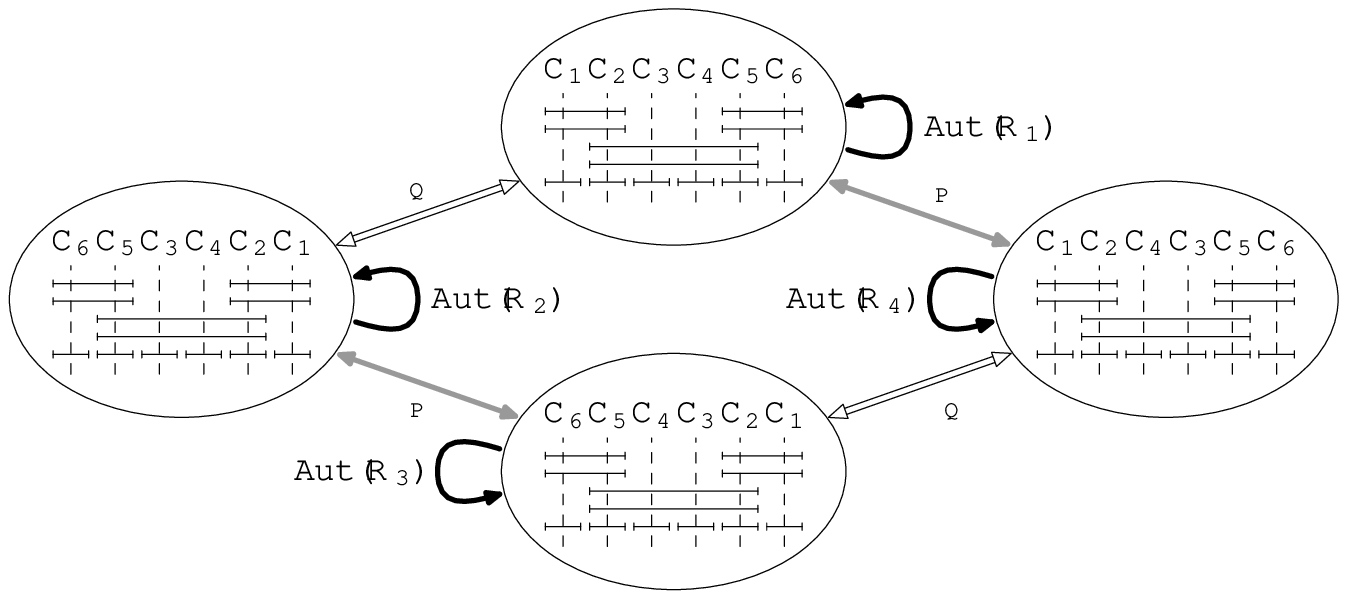}
\caption{The action of $\Aut(X)$ is transitive. An MPQ-tree $M$ of $X$ is depicted in
Fig.~\ref{fig:interval_graph_pq_tree_mpq_tree}. There are three twin classes of size two, so
$\Aut(\calR) \cong \cyc_2^3$.  The group $\Aut(M)$ is generated by $\pi_Q$ corresponding to flipping
the Q-node, and $\pi_P$ permuting the P-node. We have $\Aut(M) \cong \cyc_2^2$ and $\Aut(X) \cong
\cyc_2^3 \rtimes \cyc_2^2$.}
\label{fig:interval_action_diagram}
\end{figure}

\subsection{Direct Constructions.}

In this section, we explain Theorem~\ref{thm:aut_groups_interval}(i) by direct constructions. The
first construction answers the open problem of Hanlon~\cite{hanlon1982counting}.

\begin{lemma} \label{lem:int_to_tree}
For $X \in \int$, there exists $T \in \tree$ such that $\Aut(X) \cong \Aut(T)$.
\end{lemma}

\begin{proof}
Consider an MPQ-tree $M$ representing $X$.  We know that $\Aut(X) \cong \Ker(\Phi) \rtimes \Aut(M)$
and we inductively encode the structure of $M$ into $T$.

\begin{figure}[b!]
\centering
\includegraphics[scale=0.8]{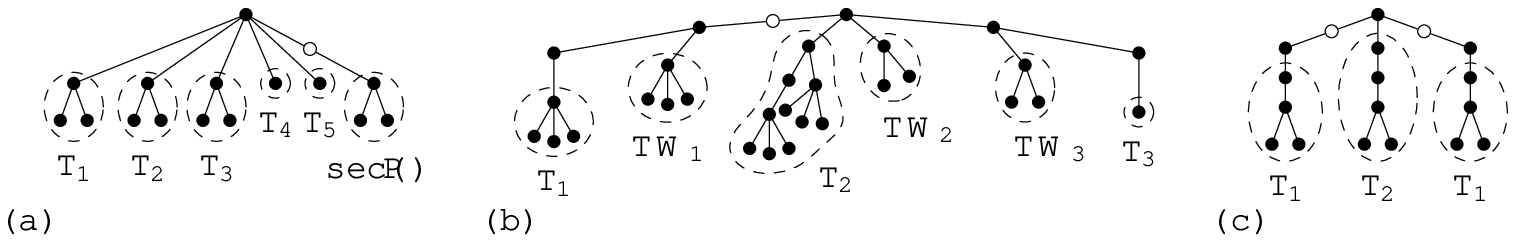}
\caption{For an interval graph $X$, a construction of a tree $T$ with $\Aut(T) \cong \Aut(X)$: (a)
The root is a P-node. (b) The root is an asymmetric Q-node. (c) The root is a symmetric Q-node.}
\label{fig:tree_construction}
\end{figure}

\emph{Case 1: The root is a P-node $P$.} Its subtrees can be encoded by trees and we just attach
them to a common root. If $\sec(P)$ is non-empty, we attach a star with $|\sec(P)|$ leaves to the
root (and we subdivide it to make it non-isomorphic to every other subtree attached to the root); see
Fig~\ref{fig:tree_construction}a. We get $\Aut(T) \cong \Aut(X)$.

\emph{Case 2: The root is a Q-node $Q$.} If $Q$ is asymmetric, we attach the trees corresponding to the
subtrees of $Q$ and stars corresponding to the vertices of twin classes in the sections of $Q$ to a
path, and possibly modify by subdivisions to make it asymmetric; see
Fig.~\ref{fig:tree_construction}b. And if $Q$ is symmetric, then $\Aut(X) \cong (G_1^2 \times G_3)
\rtimes \cyc_2$ and we just attach trees $T_1$ and $T_2$ such that $\Aut(T_i) \cong G_i$ to a path
as in Fig.~\ref{fig:tree_construction}c.  In both cases, $\Aut(T) \cong \Aut(X)$.
\end{proof}

\begin{lemma} \label{lem:tree_to_int}
For $T \in \tree$, there exists $X \in \int$ such that $\Aut(T) \cong \Aut(X)$.
\end{lemma}

\begin{proof}
For a rooted tree $T$, we construct an interval graph $X$ such that $\Aut(T) \cong \Aut(X)$ as
follows. The intervals are nested according to $T$ as shown in Fig.~\ref{fig:int_construction}. Each
interval is contained exactly in the intervals of its ancestors. If $T$ contains a vertex with only
one child, then $\Aut(T) < \Aut(X)$. This can be fixed by adding suitable asymmetric interval
graphs $Y$, as in Fig.~\ref{fig:int_construction}.
\end{proof}

\begin{figure}[t!]
\centering
\includegraphics[scale=0.8]{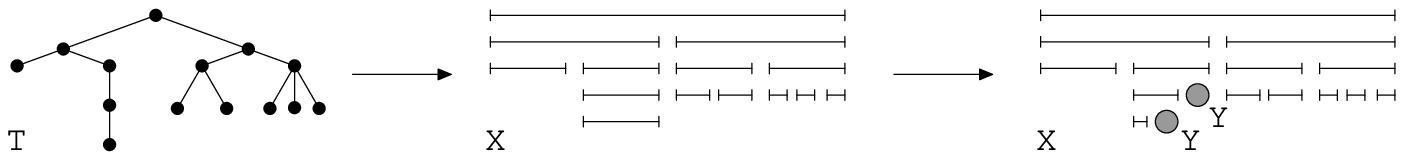}
\caption{We place the intervals following the structure of the tree. We get $\Aut(X) \cong
\sym_3 \times \sym_2 \times \sym_3$, but $\Aut(T) \cong \sym_2 \times \sym_3$.  We fix this by
attaching asymmetric interval graphs $Y$.}
\label{fig:int_construction}
\end{figure}

\subsection{Automorphism Groups of Unit Interval Graphs}

We apply the characterization of $\Aut(\int)$ derived in Theorem~\ref{thm:aut_groups_interval}(i) to
show that the automorphism groups of connected unit interval graphs are the same of caterpillars
(which form the intersection of \int\ and \tree).
The reader can make direct constructions, similarly as in Lemmas~\ref{lem:int_to_tree}
and~\ref{lem:tree_to_int}.
First, we describe $\Aut(\caterpillar)$:

\begin{lemma}\label{prop:automorphisms_caterpillars}
The class $\Aut(\caterpillar)$ consists of all groups $G_1$ and $G_1 \wr \cyc_2 \times G_2$
where $G_1$ is a direct product of symmetric groups and $G_2$ is a symmetric group.
\end{lemma}

\begin{proof}
We can easily construct caterpillars with these automorphism groups. On the
other hand, the root of an MPQ-tree $M$ representing $T$ is a Q-node $Q$ (or a P-node with at most
two children, which is trivial).  All twin classes are trivial, since $T$ is a tree. Each child of
the root is either a P-node, or a leaf. All children of a P-node are leaves.
We can determine $\Aut(X)$ as in the proof of
Theorem~\ref{thm:aut_groups_interval}(i).
\end{proof}

\begin{proof}[Proof of Theorem~\ref{thm:aut_groups_interval}(ii)]
According to Corneil \cite{corneil1995simple}, an MPQ-tree $M$ representing a connected unit interval
graph contains only one Q-node with all children as leaves. It is possible that the
sections of this Q-node are nontrivial. This equality of automorphism groups follows by
Lemma~\ref{prop:automorphisms_caterpillars} and the proof of
Theorem~\ref{thm:aut_groups_interval}(i).
\end{proof}

\section{Automorphism Groups of Circle Graphs}\label{sec:circle_graphs}

In this section, we prove Theorem \ref{thm:aut_groups_circle}. We introduce the \emph{split
decomposition} which was invented for recognizing circle graphs. We encode the split decomposition
of $X$ by a \emph{split tree} $S$ which captures all circle representations of $X$.  We define
automorphisms of $S$ and show that $\Aut(S) \cong \Aut(X)$.

\subsection{Split Decomposition}

A \emph{split} is a partition $(A,B,A',B')$ of $V(X)$ such that:
\begin{itemize}
\item For every $a \in A$ and $b \in B$, we have $ab \in E(X)$.
\item There are no edges between $A'$ and $B \cup B'$, and between $B'$ and $A \cup A'$.
\item Both sides have at least two vertices: $|A \cup A'| \ge 2$ and $|B \cup B'| \ge 2$.
\end{itemize}

The split decomposition of $X$ is constructed by taking a split of $X$ and replacing $X$ by the
graphs $X_A$ and $X_B$ defined as follows. The graph $X_A$ is created from $X[A \cup A']$ together
with a new \emph{marker vertex} $m_A$ adjacent exactly to the vertices in $A$.  The graph $X_B$ is
defined analogously for $B$, $B'$ and $m_B$; see Fig.~\ref{fig:split_graph_split-tree}a. The
decomposition is then applied recursively on $X_A$ and $X_B$. Graphs containing no splits are called
\emph{prime graphs}. We stop the split decomposition also on \emph{degenerate graphs} which are
complete graphs $K_n$ and stars $K_{1,n}$. A split decomposition is called \emph{minimal} if it is
constructed by the least number of splits. Cunningham~\cite{cunningham} proved that the minimal
split decomposition of a connected graph is unique.

The key connection between the split decomposition and circle graphs is the following: a graph $X$
is a circle graph if and only if both $X_A$ and $X_B$ are.  In a other words, a connected graph $X$
is a circle graph if and only if all prime graphs obtained by the minimal split decomposition are
circle graphs.

\begin{figure}[t]
\centering
\includegraphics[scale=0.8]{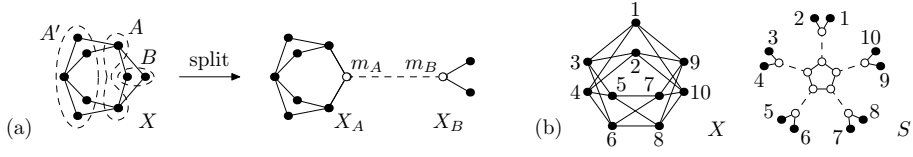}
\caption{(a) An example of a split of the graph $X$. The marker vertices are depicted in white.  The
tree edge is depicted by a dashed line.  (b) The split tree $S$ of the graph $X$. We have that
$\Aut(S) \cong \cyc_2^5 \rtimes \dih_5$.}
\label{fig:split_graph_split-tree}
\end{figure}

\heading{Split tree.} The \emph{split tree} $S$ representing a graph $X$ encodes the minimal split
decomposition. A split tree is a graph with two types of vertices (normal and marker vertices) and
two types of edges (normal and tree edges). We initially put $S = X$ and modify it according to the
minimal split decomposition. If the minimal decomposition contains a split $(A, B, A', B')$ in $Y$,
then we replace $Y$ in $S$ by the graphs $Y_A$ and $Y_B$, and connect the marker vertices $m_A$ and
$m_B$ by a \emph{tree edge} (see Fig.~\ref{fig:split_graph_split-tree}a). We repeat this recursively
on $Y_A$ and $Y_B$; see Fig.~\ref{fig:split_graph_split-tree}b. Each prime and degenerate graph is a
\emph{node} of the split tree. Since the minimal split decomposition is unique, we also have that
the split tree is unique.

Next, we prove that the split tree $S$ captures the adjacencies in $X$.

\begin{lemma}\label{lem:graph_edge_split-tree_path}
We have $xy \in E(X)$ if and only if there exists an alternating path $xm_1m_2\dots m_ky$ in $S$
such that each $m_i$ is a marker vertex and precisely the edges $m_{2i-1}m_{2i}$ are tree edges.
\end{lemma}

\begin{proof}
Suppose that $xy \in E(X)$. We prove existence of an alternating path between $x$ and $y$ by
induction according to the length of this path. If $xy \in E(S)$, then it clearly exists.  Otherwise
the split tree $S$ was constructed by applying a split decomposition. Let $Y$ be the graph in this
decomposition such that $xy \in E(Y)$ and there is a split $(A, B, A', B')$ in $Y$ in this
decomposition such that $x \in A$ and $y \in B$.  We have $x \in V(Y_A)$, $xm_A \in E(Y_A)$, $y \in
V(Y_B)$, and $ym_B \in E(Y_B)$. By induction hypothesis, there exist alternating paths between $x$
and $m_A$ and between $m_B$ and $y$ in $S$. There is a tree edge $m_A m_B$, so by joining we get an
alternating path between $x$ and $y$.  On the other hand, if there exists an alternating path $xm_1
\dots m_ky$ in $S$, by joining all splits, we get $xy \in E(X)$.
\end{proof}

\subsection{Automorphisms of Split-trees}
In~\cite{gioan2012split}, split trees are defined in terms of graph-labeled trees.  Our definition
is more suitable for automorphisms. An \emph{automorphism of a split tree $S$} is an automorphism of
$S$ which preserves the types of vertices and edges, i.e, it maps marker vertices to marker
vertices, and tree edges to tree edges. We denote the automorphism group of $S$ by $\Aut(S)$.

\begin{lemma} \label{lem:aut_split_iso_aut_circ}
If $S$ is a split tree representing a graph $X$, then $\Aut(S) \cong \Aut(X)$.
\end{lemma}

\begin{proof}
First, we show that each $\sigma \in \Aut(S)$ induces a unique automorphism $\pi$ of $X$. Since
$V(X) \subseteq V(S)$, we define $\pi = \sigma\restriction_{V(X)}$. By
Lemma~\ref{lem:graph_edge_split-tree_path}, $xy \in E(X)$ if and only if there exists an alternating
path between them in $S$. Automorphisms preserve alternating paths, so $xy \in E(X) \iff
\pi(x)\pi(y) \in E(X)$.

On the other hand, we show that $\pi \in \Aut(X)$ induces a unique automorphism $\sigma \in
\Aut(S)$. We define $\sigma \restriction_{V(X)} = \pi$ and extend it recursively on the marker
vertices. Let $(A, B, A', B')$ be a split of the minimal split decomposition in $X$.  This split is
mapped by $\pi$ to another split $(C, D, C', D')$ in the minimal split decomposition, i.e., $\pi(A)
= C$, $\pi(A') = C'$, $\pi(B) = D$, and $\pi(B') = D'$. By applying the split decomposition to the
first split, we get the graphs $X_A$ and $X_B$ with the marker vertices $m_A \in V(X_A)$ and $m_B
\in V(X_B)$. Similarly, for the second split we get $X_C$ and $X_D$ with $m_C \in V(X_C)$ and $m_D
\in V(X_D)$. Since $\pi$ is an automorphism, we have that $X_A \cong X_C$ and $X_B \cong X_D$. It
follows that the unique split trees of $X_A$ and $X_C$ are isomorphic, and similarly for $X_B$ and
$X_D$.  Therefore, we define $\sigma(m_A) = m_C$ and $\sigma(m_B) = m_D$, and we finish the rest
recursively. Since $\sigma$ is an automorphism at each step of the construction of $S$, it follows
that $\sigma \in \Aut(S)$.
\end{proof}

Similarly as for trees, there exists a \emph{center} of $S$ which is either a tree edge, or a prime
or degenerate node. If the center is a tree edge, we can modify the split tree by adding two
adjacent marker vertices in the middle of the tree edge. This clearly preserves the automorphism
group $\Aut(S)$, so from now on we assume that $S$ has a center $C$ which which is a prime or
degenerate node. We can assume that $S$ is rooted by $C$, and for a node $N$, we denote by $S[N]$
the subtree induced by $N$ and its descendants. For $N \ne C$, we call $m$ its \emph{root marker
vertex} if it is the marker vertex of $N$ attached to the parent of $N$.

\heading{Recursive Construction.}
We can describe $\Aut(S)$ recursively from the leaves to the root $C$.  Let $N$ be an arbitrary node
of $S$ and consider all its descendants.  Let $\Stab_{S[N]}(x)$ be the subgroup of $\Aut(S[N])$
which fixes $x \in V(S[N])$.  We further color the non-root marker vertices in $N$ by colors
coding isomorphism classes of the subtrees attached to them.

\begin{lemma} \label{lem:recursive_aut_split_tree}
Let $N \ne C$ be a node with the root marker vertex $m$. Let $N_1, \dots, N_k$ be the children of
$N$ with the root marker vertices $m_1,\dots,m_k$.  Then
$$\Stab_{S[N]}(m) \cong \bigl(\Stab_{S[N_1]}(m_1) \times \cdots \times \Stab_{S[N_k]}(m_k)\bigr)
\rtimes \Stab_N(m),$$
where $\Stab_N(m)$ is color preserving.
\end{lemma}

\begin{proof}
We proceed similarly as in the proof of Theorem~\ref{thm:aut_disconnected}.  We isomorphically label
the vertices of the isomorphic subtrees $S[N_i]$. Each automorphism $\pi \in \Stab_{S[N]}(m)$ is a
composition of two automorphisms $\sigma \cdot \tau$ where $\sigma$ maps each subtree $S[N_i]$ to
itself, and $\tau$ permutes the subtrees as in $\pi$ while preserving the labeling.  Therefore, the
automorphisms $\sigma$ can be identified with the elements of the direct product
$\Stab_{S[N_1]}(m_1) \times \cdots \times \Stab_{S[N_k]}(m_k)$ and the automorphisms $\tau$ with the
elements of $\Stab_N(m)$. The rest is exactly as in the proof of Theorem~\ref{thm:aut_disconnected}.
\end{proof}

The entire automorphism group $\Aut(S)$ is obtained by joining these subgroups at the central node
$C$. No vertex in $C$ has to be fixed by $\Aut(S)$.

\begin{lemma} \label{lem:recursive_aut_split_tree_root}
Let $C$ be the central node with the children $N_1,\dots,N_k$ with the root marker vertices
$m_1,\dots,m_k$. Then
$$\Aut(S) \cong \bigl(\Stab_{S[N_1]}(m_1) \times \cdots \times \Stab_{S[N_k]}(m_k)\bigr) \rtimes
\Aut(C),$$
where $\Aut(C)$ is color preserving.
\end{lemma}

\begin{proof}
Similar as the proof of Lemma~\ref{lem:recursive_aut_split_tree}.
\end{proof}

\subsection{The Action On Prime Circle Representations}\label{sec:group_acting_circle_rep}

For a circle graph $X$ with $|V(X)| = \ell$, a representation $\calR$ is completely determined by a
circular word $r_1 r_2 \cdots r_{2\ell}$ such that each $r_i \in V(X)$ and each vertex appears
exactly twice in the word.  This word describes the order of the endpoints of the chords in $\calR$
when the circle is traversed from some point counterclockwise. Two chords intersect if and only if
their occurrences alternate in the circular word. Representations are equivalent if
they have the same circular words up to rotations and reflections.

The automorphism group $\Aut(X)$ acts on the circle representations in the following way.
Let $\pi \in \Aut(X)$, then $\pi(\calR)$ is the circle representation represented by the word
$\pi(r_1) \pi(r_2) \cdots \pi(r_{2\ell})$, i.e., the chords are permuted according to $\pi$.

\begin{lemma} \label{lem:prime_circle_graphs}
Let $X$ be a prime circle graph. Then $\Aut(X)$ is isomorphic to a~subgroup of a dihedral group.
\end{lemma}

\begin{proof}
According to~\cite{gabor1989recognizing}, each prime circle graph has a unique representation
$\calR$, up to rotations and reflections of the circular order of endpoints of the chords.
Therefore, for every automorphism $\pi \in \Aut(X)$, we have $\pi(\calR) = \calR$, so $\pi$ only
rotates/reflects this circular ordering.  An automorphism $\pi \in \Aut(X)$ is called
a~\emph{rotation} if there exists $k$ such that $\pi(r_i) = r_{i+k}$, where the indexes are
used cyclically.  The automorphisms, which are not rotations, are called \emph{reflections}, since
they reverse the circular ordering.  For each reflection $\pi$, there exists $k$ such that $\pi(r_i)
= r_{k-i}$.  Notice that composition of two reflections is a rotation. Each reflection either
fixes two endpoints in the circular ordering, or none of them.

If no non-identity rotation exists, then $\Aut(X)$ is either $\cyc_1$, or $\cyc_2$.  If at least one
non-identity rotation exists, let $\rho \in \Aut(X)$ be the non-identity rotation with the smallest
value $k$, called the \emph{basic rotation}. Observe that $\langle \rho \rangle$ contains all
rotations, and if its order is at least three, then the rotations act semiregularly on $X$. If there
exists no reflection, then $\Aut(X) \cong \cyc_n$. Otherwise, $\langle \rho \rangle$ is a subgroup
of $\Aut(X)$ of index two.  Let $\varphi$ be any reflection, then $\rho \varphi \rho = \varphi$ and
$\Aut(X) \cong \dih_n$.
\end{proof}

\begin{lemma}\label{lem:prime_circle_graphs_stab} 
Let $X$ be a prime circle graph and let $m \in V(X)$. Then $\Stab_X(m)$ is isomorphic to a subgroup
of $\cyc_2^2$.
\end{lemma}

\begin{proof}
Let $mA\hat{m}B$ be a circular ordering representing $X$, where $m$ and $\hat{m}$ are the endpoints
of the chord representing $m$, and $A$ and $B$ are sequences of the endpoints of the remaining
chords. We distinguish $m$ and $\hat m$ to make the action of $\Stab_X(m)$ understandable.  Every
$\pi \in \Stab_X(m)$ either fixes both $m$ and $\hat{m}$, or swaps them.

Let $A'$ be the reflection of $A$ and $B'$ be the reflection of $B$.  If both $m$ and $\hat m$ are
fixed, then by the uniqueness this representation can only be reflected along the chord $m$. If such
an automorphism exists in $\Stab_X(m)$, we denote it by $\varphi_m$ and we have
$\varphi_m(mA\hat{m}B) = mB'\hat{m}A'$.  If $m$ and $\hat m$ are swapped, then by the uniqueness
this representation can be either reflected along the line orthogonal to the chord $m$, or by the
$180^\circ$ rotation.  If these automorphisms exist in $\Stab_X(m)$, we denote them by
$\varphi_\perp$ and $\rho$, respectively. We have $\varphi_\perp(mA\hat{m}B) = \hat{m}A'mB'$ and
$\rho(mA\hat{m}B) = \hat{m}BmA$.  Figure~\ref{fig:prime_graph_vertex_stabilizer} shows an example.

\begin{figure}[t]
\centering
\includegraphics[scale=0.8]{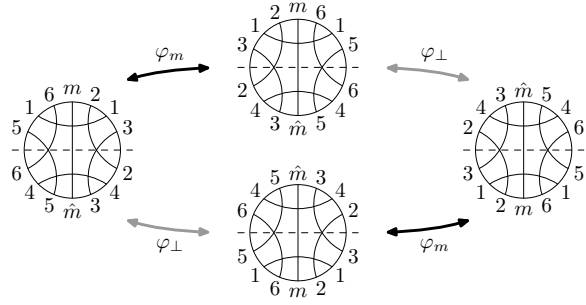}
\caption{A prime circle graph $X$ with $\Stab_X(m) \cong \cyc_2^2$.}
\label{fig:prime_graph_vertex_stabilizer}
\end{figure}

All three automorphisms $\varphi_m$, $\varphi_\perp$ and $\rho$ are involutions, and $\rho =
\varphi_\perp \cdot \varphi_m$.  Since $\Stab_X(m)$ is generated by those which exist, it is a
subgroup of $\cyc_2^2$.
\end{proof}

\subsection{The Inductive Characterization}

By Lemma~\ref{lem:aut_split_iso_aut_circ}, it is sufficient to determine the automorphism groups of
split trees. We proceed from the leaves to the root, similarly as in Theorem~\ref{thm:jordan_trees}.

\begin{lemma} \label{lem:stab_circle}
The class $\Sigma$ defined in Theorem~\ref{thm:aut_groups_circle} consists of the following groups:
\begin{equation} \label{eq:stab_circle}
\Sigma = \bigl\{G : X \in \hbox{\rm connected}\ \circle, x \in V(X), G \cong \Stab_X(x)\bigr\}.
\end{equation}
\end{lemma}

\begin{figure}[b]
\centering
\includegraphics[scale=0.8]{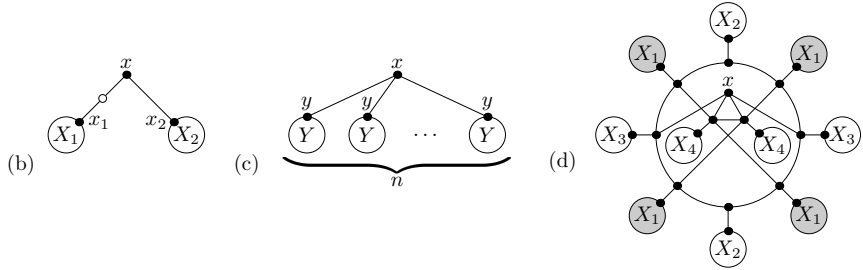}
\caption{The construction of the group in (d). The eight-cycle in $X$ can be reflected horizontally,
vertically and rotated by $180^\circ$.}
\label{fig:circle_stabilizer_groups_construction}
\end{figure}

\begin{proof}
First, we show that (\ref{eq:stab_circle}) is closed under (b) to (d); see
Fig.~\ref{fig:circle_stabilizer_groups_construction}. For (b), let $X_1$ and $X_2$ be circle graphs
such that $\Stab_{X_i}(x_i) \cong G_i$.  We construct $X$ as in
Fig.~\ref{fig:circle_stabilizer_groups_construction}b, and we get $\Stab_X(x) \cong G_1 \times G_2$.
For (c), let $Y$ be a circle graph with $\Stab_{Y}(y) \cong G$.  As $X$, we take $n$ copies of $Y$
and add a new vertex $x$ adjacent to all copies of $y$.  Clearly, we get $\Stab_X(x) \cong G \wr
\sym_n$.  For (d), let $G_1, G_2, G_3, G_4 \in \Sigma$, and let $X_i$ be a circle graph with
$\Stab_{X_i}(x_i) \cong G_i$.  We construct a graph $X$ as shown in
Fig.~\ref{fig:circle_stabilizer_groups_construction}. We get $\Stab_X(x) \cong \bigl(G_1^4 \times
G_2^2 \times G_3^2 \times G_4^2\bigr) \rtimes \cyc_2^2$.

Next we show that every group from (\ref{eq:stab_circle}) belongs to $\Sigma$.
Let $X$ be a circle graph with $x \in V(X)$, and we want to show that $\Stab_X(x) \in \Sigma$.
Since $\Aut(S) \cong \Aut(X)$ by Lemma~\ref{lem:aut_split_iso_aut_circ}, we have $\Stab_S(x) \cong
\Stab_X(x)$ where $x$ is a non-marker vertex. We prove this by induction according to the number of
nodes of $S$, for the single node it is either a subgroup $\cyc_2^2$ (by
Lemma~\ref{lem:prime_circle_graphs_stab}), or a symmetric group.

Let $N$ be the node containing $x$, we can think of it as the root and $x$ being a root marker
vertex. Therefore, by Lemma~\ref{lem:recursive_aut_split_tree}, we have
$$\Stab_{S}(x) \cong \bigl(\Stab_{S[N_1]}(m_1) \times \cdots \times \Stab_{S[N_k]}(m_k)\bigr) \rtimes
\Stab_N(x),$$
where $N_1,\dots,N_k$ are the children of $N$ and $m_1,\dots,m_k$ their root marker vertices.
By the induction hypothesis, $\Stab_{S[N_i]}(m_i) \in \Sigma$. There are two cases:

\emph{Case 1: $N$ is a degenerate node.}
Then $\Stab_N(x)$ is a direct product of symmetric groups. The subtrees attached to marker vertices
of each color class can be arbitrarily permuted, independently of each other. Therefore
$\Stab_{S}(x)$ can be constructed using (b) and (c), exactly as in
Theorem~\ref{thm:aut_disconnected}.

\emph{Case 2: $N$ is a prime node.}
By Lemma~\ref{lem:prime_circle_graphs_stab}, $\Stab_{N}(x)$ is a subgroup of $\cyc_2^2$.
When it is trivial or $\cyc_2$, observe that $\Stab_{S}(x)$ can be constructed using (b) and (c).
The only remaining case is when it is $\cyc_2^2$. The action of $\cyc_2^2$ on $V(N)$ can have
orbits of sizes $4$, $2$, and $1$. By Orbit-Stabilizer Theorem, each orbit of size $2$
has also a stabilizer of size $2$, having exactly one non-trivial element. Therefore, there are at
most three types of orbits of size $2$, according to which of elements $(1,0)$, $(0,1)$ and $(1,1)$
stabilizes them. Figure~\ref{fig:circle_stabilizer_groups_construction} shows that all three types
of orbits are possible.

Let $G_1$ be the direct product of all $\Stab_{S[N_i]}(m_i)$, one from each orbit of size four.
The groups $G_2$, $G_3$, and $G_4$ are defined similarly for the three types of orbits of size two,
and $G_5$ for the orbits of size one. We get that
$$\Stab_{S}(x) \cong \bigl (G_1^4 \times G_2^2 \times G_3^2 \times G_4^2\bigr)
\rtimes_{\varphi} \cyc_2^2 \times G_5,$$
where $\varphi(1,0)$ and $\varphi(0,1)$ swap the coordinates as the horizontal and vertical
reflections in Fig.~\ref{fig:circle_stabilizer_groups_construction}d, respectively.  Thus
$\Stab_{S}(x)$ can be build using (b) and (d).
\end{proof}

Now, we prove Theorem~\ref{thm:aut_groups_circle}.

\begin{figure}[b]
\centering
\includegraphics[scale=0.8]{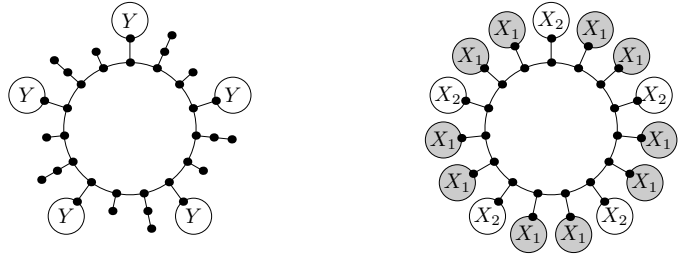}
\caption{The construction of the described groups.}
\label{fig:circle_groups_construction}
\end{figure}

\begin{proof}[Proof of Theorem~\ref{thm:aut_groups_circle}]
We first prove that $\Autcc$ contains all described groups. Let $G \in \Sigma$ and let $Y$ be a
connected circle graph with $\Stab_{Y}(y) \cong G$.  We take $n$ copies of $Y$ and attach them by
$y$ to the graph depicted in Fig.~\ref{fig:circle_groups_construction} on the left.  Clearly, we get
$\Aut(x) \cong G^n \rtimes \cyc_n$.  Let $G_1, G_2 \in \Sigma$ and let $X_1$ and $X_2$ be connected
circle graphs such that $\Stab_{X_i}(x_i) \cong G_i$ and $X_1 \not\cong X_2$. We construct
a graph $X$ by attaching $n$ copies of $X_1$ by $x_1$ and $2n$ copies of $X_2$ by $x_2$ as in
Fig.~\ref{fig:circle_groups_construction} on the right. We get $\Aut(X) \cong (G_1^n \times
G_2^{2n}) \rtimes \dih_n$.

Let $X$ be a connected circle graph, we want to show that $\Aut(X)$ can be constructed in the above
way. Let $S$ be its split, by Lemma~\ref{lem:aut_split_iso_aut_circ} we have $\Aut(S) \cong
\Aut(X)$. For the central node $C$, we get by Lemma~\ref{lem:recursive_aut_split_tree_root} that
$$\Aut(S) \cong \bigl(\Stab_{S[N_1]}(m_1) \times \cdots \times \Stab_{S[N_k]}(m_k)\bigr) \rtimes
\Aut(C),$$
where $N_1,\dots,N_k$ are children of $C$ and $m_1,\dots,m_k$ are their root marker vertices.
By Lemma~\ref{lem:stab_circle}, we know that each $\Stab_{S[N_i]} \in \Sigma$ and also $\prod
\Stab_{S[N_i]}(m_i) \in \Sigma$. The rest follows by analysing the automorphism group $\Aut(C)$ and
its orbits.

\emph{Case 1: $C$ is a degenerate node.} This is exactly the same as Case 1 in the proof of
Lemma~\ref{lem:stab_circle}. We get that $\Aut(S) \in \Sigma$, so it is the semidirect product with
$\cyc_1$.

\emph{Case 2: $C$ is a prime node.} By Lemma~\ref{lem:prime_circle_graphs}, we know that $\Aut(C)$
is isomorphic to either $\cyc_n$, or $\dih_n$. If $n \le 2$, we can show by a similar argument that
$\Aut(S) \in \Sigma$.

If $\Aut(C) \cong \cyc_n$, where $n \ge 3$, then by Lemma~\ref{lem:prime_circle_graphs} we know that
$\Aut(C)$ consists of rotations which act semiregularly. Therefore each orbit of $\Aut(C)$ is of
size $n$ and $\Aut(C)$ acts isomorphically on them. Let $G \in \Sigma$ be the direct product of
$\Stab_{S[N_i]}(m_i)$, one for each orbit of $\Aut(C)$. It follows that $$\Aut(S) \cong G^n \rtimes
\Aut(C) = G \wr \cyc_n.$$

If $\Aut(C) \cong \dih_n$, where $n \ge 3$, then by Lemma~\ref{lem:prime_circle_graphs} there exists
a subgroup of rotations of index two, acting semiregularly. Therefore each orbit of $\Aut(C)$ is of
size $n$ or $2n$. On the orbits of size $2n$, we know that $\Aut(C)$ acts regularly. Let $\rho \in
\Aut(C)$ be the basic rotation by $k$. Then the chords belonging to an orbit of size $n$ are
cyclically shifted by $k$ endpoints. Therefore $\Aut(C)$ acts on all of them isomorphically, exactly
as on the vertices of a regular $n$-gon. Let $G_1 \in \Sigma$ be the direct product of
$\Stab_{S[N_i]}(m_i)$, one from each orbit of size $n$, and let $G_2 \in \Sigma$ be the direct
product of $\Stab_{S[N_i]}(m_i)$, one for each orbit of size $2n$. We get:
$$\Aut(S) \cong (G_1^n \times G_2^{2n}) \rtimes \Aut(C) = (G_1^n \times G_2^{2n}) \rtimes \dih_n,$$
where $\dih_n$ permutes the coordinates in $G_1^n$ exactly as the vertices of a regular $n$-gon,
and permutes the coordinates in $G_2^{2n}$ regularly.
\end{proof}

\subsection{The Action on Circle Representations.}

For a connected circle graph $X$, the set $\Repfac$ consists of all circular orderings of the
endpoints of the chords which give a correct representation of $X$. Then $\pi(\calR)$ is the
representation in which the endpoints are mapped by $\pi$. The stabilizer $\Aut(\calR)$ can only
rotate/reflect this circular ordering, so it is a subgroup of a dihedral group. For prime
circle graphs, we know that $\Aut(\calR) = \Aut(X)$. A general circle graph may have many different
representations, and the action of $\Aut(X)$ on them may consist of several non-isomorphic orbits and
$\Aut(\calR)$ may not be a normal subgroup of $\Aut(X)$. 

The above results have the following interpretation in terms of the action of $\Aut(X)$. By
Lemma~\ref{lem:aut_split_iso_aut_circ}, we know that $\Aut(S) \cong \Aut(X)$. We assume that the
center $C$ is a prime circle graph, otherwise $\Aut(\calR)$ is very restricted ($\cyc_1$ or
$\cyc_2$) and not very interesting.  We choose a representation $\calR$ belonging to the smallest
orbit, i.e., $\calR$ is one of the most symmetrical representations.  Then $\Aut(\calR)$ consists of 
the rotations/reflections of $C$ described in the proof of Theorem~\ref{thm:aut_groups_circle}. 

The action of $\Aut(X)$ on this orbit is described by the point-wise stabilizer $H$ of $C$ in
$\Aut(S)$. We know that $H = \prod \Stab_{S[N_i]}(m_i)$ as described in
Lemma~\ref{lem:stab_circle}. When $N_i$ is a prime graph, we can apply reflections and rotations
described in Lemma~\ref{lem:prime_circle_graphs_stab}, so we get a subgroup of $\cyc_2^2$.  If $N_i$
is a degenerate graph, then isomorphic subtrees can be arbitrarily permuted which corresponds to
permuting small identical parts of a circle representation.
It follows that $\Aut(X) \cong H \rtimes \Aut(\calR)$.

\section{Automorphism Groups of Comparability and Permutation Graphs} \label{sec:comparability_graphs}

All transitive orientations of a graph are efficiently captured by the modular decomposition which
we encode into the modular tree.  We study the induced action of $\Aut(X)$ on the set of all
transitive orientations.  We show that this action is captured by the modular tree, but for general
comparability graphs its stabilizers can be arbitrary groups.  In the case of permutation graphs, we
study the action of $\Aut(X)$ on the pairs of orientations of the graph and its complement, and show
that it is semiregular. Using this, we prove Theorem~\ref{thm:aut_groups_perm}. We also show that an
arbitrary graph can be encoded into a comparability graph of the dimension at most four, which
establishes Theorem~\ref{thm:kdim_aut_groups_and_gi}.

\subsection{Modular Decomposition} \label{sec:modular_decomposition}

A \emph{module} $M$ of a graph $X$ is a set of vertices such that each $x \in V(X) \setminus M$ is
either adjacent to all vertices in $M$, or to none of them.  Modules generalize connected
components, but one module can be a proper subset of another one.  Therefore, modules lead to a
recursive decomposition of a graph, instead of just a partition. See Fig.~\ref{fig:modules}a for
examples. A module $M$ is called \emph{trivial} if $M=V(X)$ or $|M|=1$, and \emph{non-trivial}
otherwise.

If $M$ and $M'$ are two disjoint modules, then either the edges between $M$ and $M'$ form the
complete bipartite graph, or there are no edges at all; see Fig.~\ref{fig:modules}a.  In the former
case, $M$ and $M'$ are called \emph{adjacent}, otherwise they are \emph{non-adjacent}.

\begin{figure}[b]
\centering
\includegraphics[scale=0.8]{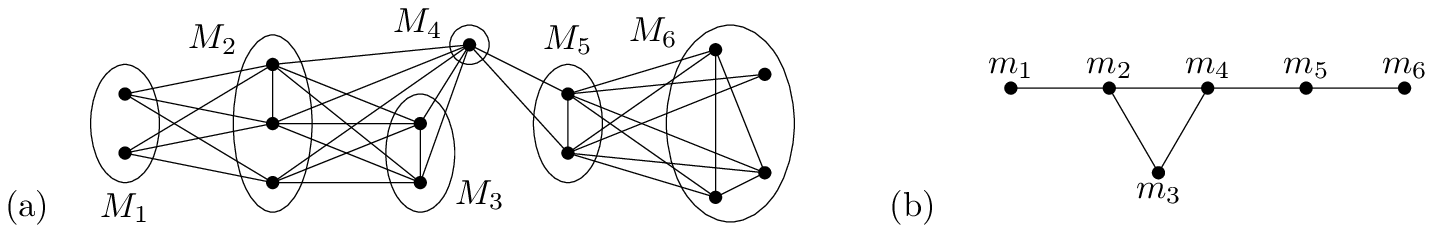}
\caption{(a) A graph $X$ with a modular partition $\calP$. (b) The quotient graph $X/\calP$ is prime.}
\label{fig:modules}
\end{figure}

\heading{Quotient Graphs.}
Let $\calP = \{M_1, \dots, M_k\}$ be a \emph{modular partition} of $V(X)$, i.e., each $M_i$ is a
module of $X$, $M_i \cap M_j = \emptyset$ for every $i\neq j$, and $M_1 \cup \cdots \cup M_k =
V(X)$. We define the \emph{quotient graph} $X/\calP$ with the vertices $m_1,\dots,m_k$ corresponding
to $M_1,\dots,M_k$ where $m_im_j \in E(X/\calP)$ if and only if $M_i$ and $M_j$ are adjacent.  In
other words, the quotient graph is obtained by contracting each module $M_i$ into the single vertex
$m_i$; see Fig.~\ref{fig:modules}b.

\heading{Modular Decomposition.}
To decompose $X$, we find some modular partition $\calP = \{M_1, \dots, M_k\}$, compute $X / \calP$
and recursively decompose $X / \calP$ and each $X[M_i]$. The recursive process terminates on
\emph{prime graphs} which are graphs containing only trivial modules.  There might be many such
decompositions for different choices of $\calP$ in each step.  In 1960s,
Gallai~\cite{gallai1967transitiv} described the \emph{modular decomposition} in which special
modular partitions are chosen and which encodes all other decompositions.

The key is the following observation. Let $M$ be a module of $X$ and let $M' \subseteq M$. Then $M'$
is a module of $X$ if and only if it is a module of $X[M]$.  A graph $X$ is called \emph{degenerate}
if it is $K_n$ or $\overline{K}_n$.	We construct the modular decomposition of a graph $X$ in the
following way, see Fig.~\ref{fig:modular_tree}a for an example:
\begin{itemize}
\item If $X$ is a prime or a degenerate graph, then we terminate the modular decomposition on $X$.
We stop on degenerate graphs since every subset of vertices forms a module, so it is not useful to
further decompose them.
\item Let $X$ and $\overline{X}$ be connected graphs.  Gallai~\cite{gallai1967transitiv} shows that
the inclusion maximal proper subsets of $V(X)$ which are modules form a modular partition $\calP$ of
$V(X)$, and the quotient graph $X/\calP$ is a prime graph; see Fig.~\ref{fig:modules}. We
recursively decompose $X[M]$ for each $M \in \calP$.
\item If $X$ is disconnected and $\overline{X}$ is connected, then every union of connected
components is a module. Therefore the connected components form a modular partition $\calP$ of
$V(X)$, and the quotient graph $X/\calP$ is an independent set. We recursively decompose $X[M]$ for
each $M \in \calP$.
\item If $\overline{X}$ is disconnected and $X$ is connected, then the modular decomposition is
defined in the same way on the connected components of $\overline{X}$.  They form a modular
partition $\calP$ and the quotient graph $X/\calP$ is a complete graph. We recursively decompose
$X[M]$ for each $M \in \calP$.
\end{itemize}

\subsection{Modular Tree.}

We encode the modular decomposition by the \emph{modular tree} $T$, similarly as
the split decomposition is captured by the split tree in Section~\ref{sec:circle_graphs}.  The
modular tree $T$ is a graph with two types of vertices (normal and \emph{marker vertices}) and two
types of edges (normal and \emph{directed tree edges}). The directed tree edges connect the prime
and degenerate graphs encountered in the modular decomposition (as quotients and terminal graphs)
into a rooted tree. 

We give a recursive definition.  Every modular tree has an induced subgraph called
\emph{root node}.   If $X$ is a prime or a degenerate graph, we define $T = X$ and its root node is
equal $T$.  Otherwise, let $\calP = \{M_1,\dots,M_k\}$ be the used modular partition of $X$ and let
$T_1,\dots,T_k$ be the corresponding modular trees for $X[M_1],\dots,X[M_k]$. The modular tree $T$
is the disjoint union of $T_1,\dots,T_k$ and of the quotient $X / \calP$ with the
marker vertices $m_1,\dots,m_k$. To every graph $T_i$, we add a new marker vertex $m'_i$ such that
$m'_i$ is adjacent exactly to the vertices of the root node of $T_i$. We further add a tree edge
oriented from $m_i$ to $m'_i$.  For an example, see Fig.~\ref{fig:modular_tree}b.

\begin{figure}[t]
\centering
\includegraphics[scale=0.8]{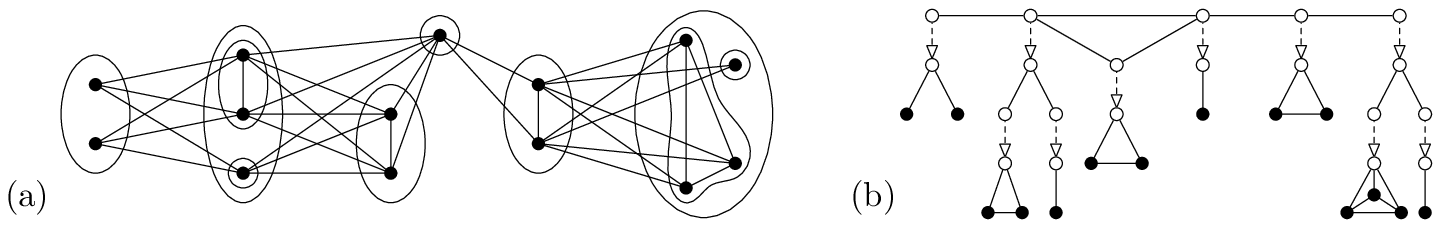}
\caption{(a) The graph $X$ from Fig.~\ref{fig:modules} with the modular partitions used in the
modular decomposition. (b) The modular tree $T$ of $X$, the marker vertices are white, the tree
edges are dashed.}
\label{fig:modular_tree}
\end{figure}

The modular tree of $X$ is unique.  The graphs encountered in the modular decomposition are called
\emph{nodes of $T$}, or alternatively root nodes of some modular tree in the construction of $T$.
For a node $N$, its subtree is the modular tree which has $N$ as the root node. \emph{Leaf nodes}
correspond to the terminal graphs in the modular decomposition, and \emph{inner nodes} are the
quotients in the modular decomposition. All vertices of $X$ are in leaf nodes and all marker
vertices, corresponding to modules of $X$, are in inner nodes.

Similarly as in Lemma~\ref{lem:graph_edge_split-tree_path}, the modular tree $T$ captures the
adjacencies in $X$.

\begin{lemma}\label{lem:adjacency_modular_tree}
We have $xy \in E(X)$ if and only if there exists an alternating path $xm_1m_2\dots m_ky$ in the
modular tree $T$ such that each $m_i$ is a marker vertex and precisely the edges $m_{2i-1}m_{2i}$
are tree edges.
\end{lemma}

\begin{proof}
Both $x$ and $y$ belong to leaf nodes. If there exists an alternating path, let $N$ be the node
which is the common ancestor of $x$ and $y$. This path has an edge $m_{2i}m_{2i+1}$ in $N$.
These vertices correspond to adjacent modules $M_{2i}$ and $M_{2i+1}$ such that $x \in M_{2i}$ and
$y \in M_{2i+1}$. Therefore $xy \in E(X)$.

On the other hand, let $N$ be the common ancestor of $x$ and $y$, such that $m_x$ is the marker
vertex on a path from $x$ to $N$ and similarly $m_y$ is the marker vertex for $y$ and $N$. If $xy
\in E(X)$, then the corresponding modules $M_x$ and $M_y$ has to be adjacent, so we can construct an
alternating path from $x$ to $y$.
\end{proof}

\subsection{Automorphisms of Modular Trees.}

An automorphism of the modular tree $T$ has to preserve the types of vertices and edges and the
orientation of tree edges. We denote the automorphism group of $T$ by $\Aut(T)$.

\begin{lemma}\label{lem:aut_group_modular_tree}
If $T$ is the modular tree of a graph $X$, then $\Aut(X) \cong \Aut(T)$.
\end{lemma}

\begin{proof}
First, we show that each automorphism $\sigma \in \Aut(T)$ induces a unique automorphism of $X$.
Since $V(X) \subseteq V(T)$, we define $\pi = \sigma\restriction_{V(X)}$. By
Lemma~\ref{lem:adjacency_modular_tree}, $xy \in E(X)$ if and only if there exists an alternating
path in $T$ connecting them. Automorphisms preserve alternating paths, so $xy \in E(X) \iff
\pi(x)\pi(y) \in E(X)$.

For the converse, we prove that $\pi \in \Aut(X)$ induces a unique automorphism $\sigma \in
\Aut(T)$. We define $\sigma \restriction_{V(X)} = \pi$ and extend it recursively on the marker
vertices. Let $\calP = \{M_1, \dots, M_k\}$ be the modular partition of $X$ used in the modular
decomposition. It is easy to see that $\Aut(X)$ induces an action on $\calP$. If $\pi(M_i)
= M_j$, then clearly $X[M_i]$ and $X[M_j]$ are isomorphic. We define $\sigma(m_i) = m_j$ and
$\sigma(m_i') = m_j'$, and finish the rest recursively. Since $\sigma$ is an automorphism at each
step of the construction, it follows that $\sigma \in \Aut(T)$.
\end{proof}

\heading{Recursive Construction.}
We can build $\Aut(T)$ recursively. Let $N$ be the root node of $T$.  Suppose that we know the
automorphism groups $\Aut(T_1),\dots,\Aut(T_k)$ of the subtrees $T_1,\dots,T_k$ of all children of
$N$.  We further color the marker vertices in $N$ by colors coding isomorphism classes of the
subtrees $T_1,\dots,T_k$.

\begin{lemma} \label{lem:recursive_aut_modular_tree}
Let $N$ be the root node of $T$ with subtrees $T_1,\dots,T_k$. Then
$$\Aut(T) \cong \bigl(\Aut(T_1) \times \cdots \times \Aut(T_k)\bigr) \rtimes \Aut(N),$$
where $\Aut(N)$ is color preserving.
\end{lemma}

\begin{proof}
Recall the proof of Theorem~\ref{thm:aut_disconnected}.  We isomorphically label
the vertices of the isomorphic subtrees $T_i$. Each automorphism $\pi \in \Aut(T)$ is a composition
of two automorphisms $\sigma \cdot \tau$ where $\sigma$ maps each subtree $T_i$ to itself, and
$\tau$ permutes the subtrees as in $\pi$ while preserving the labeling.  Therefore, the
automorphisms $\sigma$ can be identified with the elements of 
$\Aut(T_1) \times \cdots \times \Aut(T_k)$ and the automorphisms $\tau$ with the elements of
$\Aut(N)$.  The rest is exactly as in the proof of Theorem~\ref{thm:aut_disconnected}.
\end{proof}

With no further assumptions on $X$, if $N$ is a prime graph, then $\Aut(N)$ can be isomorphic to an
arbitrary group, as shown in Section~\ref{sec:kdim}. If $N$ is a degenerate graph, then $\Aut(N)$ is
a direct product of symmetric groups.

\heading{Automorphism Groups of Interval Graphs.} In Section~\ref{sec:interval_graphs}, we proved
using MPQ-trees that $\Aut(\int) = \Aut(\tree)$. The modular decomposition gives an alternative
derivation that $\Aut(\int) \subseteq \Aut(\tree)$ by Lemma~\ref{lem:recursive_aut_modular_tree} and
the following:

\begin{lemma}
For a prime interval graph $X$, $\Aut(X)$ is a subgroup of $\cyc_2$.
\end{lemma}

\begin{proof}
Hsu~\cite{hsu1995m} proved that prime interval graphs have exactly two consecutive orderings of the
maximal cliques. Since $X$ has no twin vertices, $\Aut(X)$ acts semiregularly on the consecutive
orderings and there is at most one non-trivial automorphism in $\Aut(X)$.
\end{proof}

\subsection{Automorphism Groups of Comparability Graphs} \label{sec:aut_comp}

In this section, we explain the structure of the automorphism groups of comparability graphs, in
terms of actions on sets of transitive orientations.

\heading{Structure of Transitive Orientations.}
Let $\to$ be a transitive orientation of $X$ and let $T$ be the modular tree. For
modules $M_1$ and $M_2$, we write $M_1 \to M_2$ if $x_1 \to x_2$ for all $x_1 \in M_1$ and $x_2 \in
M_2$.  Gallai~\cite{gallai1967transitiv} shows the following properties. If $M_1$ and $M_2$ are adjacent
modules of a partition used in the modular decomposition, then either $M_1 \to M_2$, or $M_1
\leftarrow M_2$. The graph $X$ is a comparability graph if and only if each node
of $T$ is a comparability graph. Every prime comparability graph has exactly two transitive
orientations, one being the reversal of the other.

The modular tree $T$ encodes all transitive orientations as follows. For each prime node of $T$, we
arbitrarily choose one of the two possible orientations. For each degenerate node, we choose some
orientation.  (Where $K_n$ has $n!$ possible orientations and $\overline{K_n}$ has the unique
orientation.) A transitive orientation of $X$ is then constructed as follows.  We orient the
edges of leaf nodes as above. For a node $N$ partitioned in the modular decomposition by $\calP =
\{M_1,\dots,M_k\}$, we orient $X[M_i] \to X[M_j]$ if and only if $m_i \to m_j$ in $N$.  It is easy
to check that this gives a valid transitive orientation, and every transitive orientation can be
constructed by some orientation of the nodes of $T$. We note that this implies that the dimension of
the transitive orientation is the maximum of the dimensions over all nodes of $T$, and that this
dimension is the same for every transitive orientation.

\heading{Action Induced On Transitive Orientations.} Let $\orient(X)$ be the set
of all transitive orientations of $X$.  Let $\pi \in \Aut(X)$ and $\toord \in
\orient(X)$. We define the orientation $\pi(\toord)$ as follows:
$$x \to y \implies \pi(x) \mathrel{\pi(\toord)} \pi(y),\qquad \forall x,y \in V(X).$$
We can observe that $\pi(\toord)$ is a transitive orientation of $X$, so $\pi(\toord) \in
\orient(X)$; see Fig.~\ref{fig:comparability_modular_tree}. It easily follows that $\Aut(X)$
defines an action on $\orient(X)$.

\begin{figure}[b]
\centering
\includegraphics[scale=0.8]{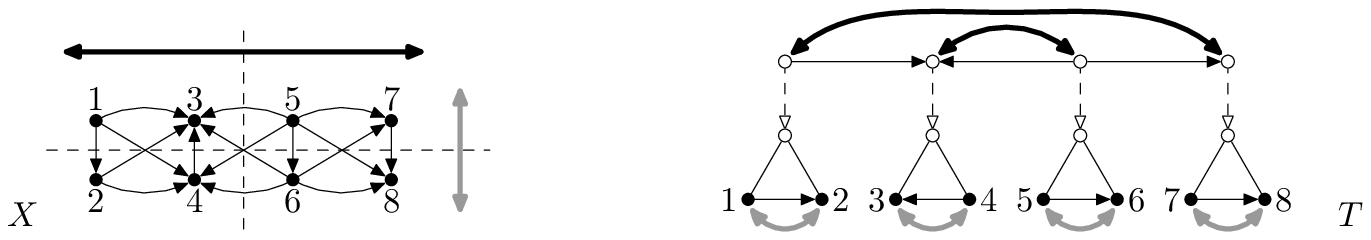}
\caption{Two automorphisms reflect $X$ and change the transitive orientation, and their action on
the modular tree $T$.}
\label{fig:comparability_modular_tree}
\end{figure}

Let $\Stab(\toord)$ be the stabilizer of some orientation $\toord \in \orient(X)$. It consists of
all automorphisms which preserve this orientation, so only the vertices that are incomparable in
$\toord$ can be permuted. In other words, $\Stab(\toord)$ is the automorphism group of the poset
created from the transitive orientation $\toord$ of $X$. Since posets are
universal~\cite{universal_posets,universal_posets2}, $\Stab(\toord)$ can be arbitrary groups and in
general the structure of $\Aut(X)$ cannot be derived from its action on $\orient(X)$, which is not
faithful enough.

Lemma~\ref{lem:recursive_aut_modular_tree} allows to understand it in terms of $\Aut(T)$ for the
modular tree $T$ representing $X$.  Each automorphism of $\Aut(T)$ somehow acts inside each node, and
somehow permutes the attached subtrees. Consider a node $N$ with attached subtrees $T_1,\dots,T_k$.
If $\sigma \in \Stab(\toord)$, then it preserves the orientation in $N$. Therefore if it maps $T_i$
to $\sigma(T_i)$, the corresponding marker vertices are necessarily incomparable in $N$.  If $N$ is
an independent set, the isomorphic subtrees can be arbitrarily permuted in $\Stab(\toord)$. If $N$
is a complete graph, all subtrees are preserved in $\Stab(\toord)$.  If $N$ is a prime graph, then
isomorphic subtrees of incomparable marker vertices can be permuted according to the structure of
$N$ which can be complex.

It is easy to observe that stabilizers of all orientations are the same and that $\Stab(\toord)$
is a normal subgroup. Let $H = \Aut(X) / \Stab(\toord)$, so $H$ captures the action of
$\Aut(X)$ on $\orient(X)$. This quotient group can be constructed recursively from the structure of
$T$, similarly to Lemma~\ref{lem:recursive_aut_modular_tree}.  Suppose that we know $H_1,\dots,H_k$
of the subtrees $T_1,\dots,T_k$. If $N$ is an independent set, there is exactly one transitive
orientation, so $H \cong H_1 \times \cdots \times H_k$. If $N$ is a complete graph, isomorphic
subtrees can be arbitrarily permuted, so $H$ can be constructed exactly as in
Theorem~\ref{thm:aut_disconnected}.  If $N$ is a prime node, there are exactly two transitive
orientations. If there exists an automorphism changing the orientation of $N$, we can describe $H$
by a semidirect product with $\cyc_2$ as in Theorem~\ref{thm:aut_disconnected}. And if $N$ is
asymmetric, then $H \cong H_1 \times \cdots \times H_k$. In particular, this description implies that
$H \in \Aut(\tree)$.

\subsection{Automorphism Groups of Permutation Graphs}

In this section, we derive the characterization of $\Aut(\perm)$ stated in Theorem~\ref{thm:aut_groups_perm}.

\heading{Action Induced On Pairs of Transitive Orientations.} Let $X$ be a permutation graph. In
comparison to general comparability graphs, the main difference is that both $X$ and $\overline{X}$
are comparability graphs. From the results of Section~\ref{sec:aut_comp} it follows that
$\Aut(X)$ induces an action on both $\orient(X)$ and $\orient(\overline{X})$.  Let
$\orient(X,\overline{X}) = \orient(X) \times \orient(\overline{X})$, and we work with one action on
the pairs $(\toord,\overtoord) \in \orient(X,\overline{X})$.
Figure~\ref{fig:action_orientations} shows an example.

\begin{figure}[b]
\centering
\includegraphics[scale=0.8]{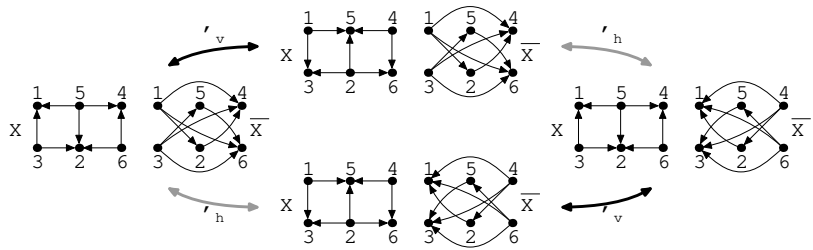}
\caption{The action of $\Aut(X)$ on four pairs of transitive orientations $X$. The black generator
flips the orientation of $X$, the gray automorphism of both $X$ and $\overline X$.}
\label{fig:action_orientations}
\end{figure}

\begin{lemma}\label{lem:action_semiregular}
For a permutation graph $X$, the action of $\Aut(X)$ on $\orient(X,\overline X)$ is semiregular.
\end{lemma}

\begin{proof}
Since a permutation belonging to the stabilizer of $(\toord,\overtoord)$ fixes both
orientations, it can only permute incomparable elements. But incomparable elements in $\toord$ are
exactly the comparable elements in $\overtoord$, so the stabilizer is trivial.
\end{proof}

\begin{lemma}\label{lem:groups_prime_perm}
For a prime permutation graph $X$, $\Aut(X)$ is a subgroup of $\cyc_2^2$.
\end{lemma}

\begin{proof}
There are at most four pairs of orientations in $\orient(X,\overline X)$, so by
Lemma~\ref{lem:action_semiregular} the order of $\Aut(X)$ is at most four. If $\pi \in \Aut(X)$,
then $\pi^2$ fixes the orientations of both $X$ and $\overline X$. Therefore $\pi^2$ belongs to the
stabilizers and it is an identity. Thus $\pi$ is the involution and $\Aut(X)$ is a subgroup of
$\cyc_2^2$.
\end{proof}

\begin{figure}[t]
\centering
\includegraphics[scale=0.8]{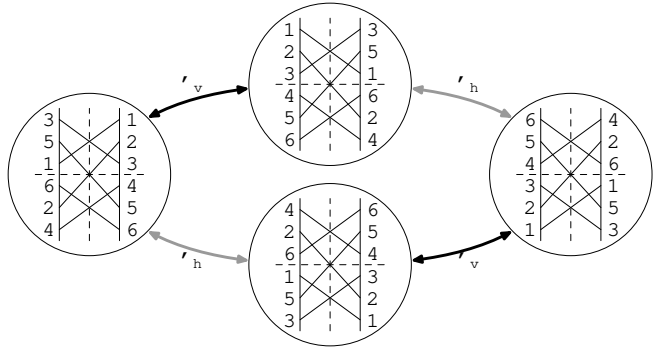}
\caption{Four representations of a symmetric permutation graph. The black automorphism is the
horizontal reflection, the gray automorphism is the vertical reflection.}
\label{fig:action_representations}
\end{figure}

\heading{Geometric Interpretation.} First, we explain the result $\perm =
\dim2$ of Even et al.~\cite{even1972permutation}.  Let $\toord \in \orient(X)$ and
$\overtoord \in \orient(\overline X)$, and let $\overtoord_R$ be the reversal of
$\overtoord$. We construct two linear orderings $L_1 = \toord \cup \overtoord$ and
$L_2 = \toord \cup \overtoord_R$. The comparable pairs in $L_1 \cap L_2$ are precisely the
edges $E(X)$.

Consider a permutation representation of a symmetric prime permutation graph.  The vertical
reflection $\varphi_v$ corresponds to exchanging $L_1$ and $L_2$, which is equivalent to reversing
$\overtoord$. The horizontal reflection $\varphi_h$ corresponds to reversing both $L_1$ and $L_2$,
which is equivalent to reversing both $\toord$ and $\overtoord$. We denote the central $180^\circ$
rotation by $\rho = \varphi_h \cdot \varphi_v$ which corresponds to reversing $\toord$; see
Fig.~\ref{fig:action_representations}. 

\heading{The Inductive Characterization.} Now, we are ready to prove Theorem~\ref{thm:aut_groups_perm}.

\begin{proof}[Proof of Theorem~\ref{thm:aut_groups_perm}]\label{pf:aut_groups_perm}
First, we show that $\Aut(\perm)$ is closed under (b) to (d). For (b), let $G_1, G_2 \in
\Aut(\perm)$, and let $X_1$ and $X_2$ be two permutation graphs such that $\Aut(X_i) \cong G_i$. We
construct $X$ by attaching $X_1$ and $X_2$ as in Fig.~\ref{fig:counterexample}b. Clearly, $\Aut(X)
\cong G_1 \times G_2$.  For (c), let $G \in \Aut(\perm)$ and let $Y$ be a connected permutation
graph such that $\Aut(Y) \cong G$. We construct $X$ as the disjoint union of $n$ copies of $Y$; see
Fig.~\ref{fig:counterexample}c. We get $\Aut(X) \cong G \wr \sym_n$. Let $G_1, G_2, G_3 \in
\Aut(\perm)$, and let $X_1$, $X_2$, and $X_3$ be permutation graphs such that $\Aut(X_i) \cong G_i$.
We construct $X$ as in Fig.~\ref{fig:counterexample}d. We get $\Aut(X) \cong \bigl (G_1^4
\times G_2^2 \times G_3^2 \bigr) \rtimes \cyc_2^2$.

\begin{figure}[b]
\centering
\includegraphics[scale=0.8]{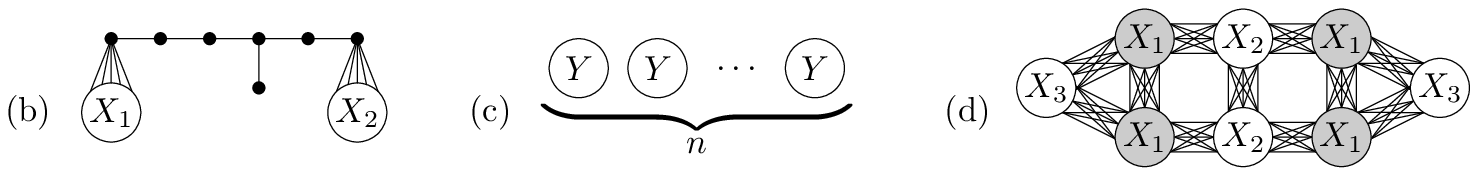}
\caption{The constructions in the proof of Theorem~\ref{thm:aut_groups_perm}.}
\label{fig:counterexample}
\end{figure}

We show the other implication by induction.  Let $X$ be a permutation graph and let $T$ be the
modular tree representing $X$. By Lemma~\ref{lem:aut_group_modular_tree}, we know that $\Aut(T)
\cong \Aut(X)$. Let $N$ be the root node of $T$, and let $T_1, \dots, T_k$ be the subtrees
attached to $N$. By the induction hypothesis, we assume that $\Aut(T_i) \in \Aut(\perm)$.  By
Lemma~\ref{lem:recursive_aut_modular_tree},
$$\Aut(T) \cong \bigl(\Aut(T_1) \times \cdots \times \Aut(T_k)\bigr) \rtimes \Aut(N).$$

\emph{Case 1: $N$ is a degenerate node.} Then $\Aut(N)$ is a direct product of symmetric groups. The
subtrees attached to marker vertices of each color class can be arbitrarily permuted, independently
of each other. Therefore $\Aut(T)$ can be constructed using (b) and (c), exactly as in
Theorem~\ref{thm:aut_disconnected}.

\emph{Case 2: $N$ is a prime node.} By Lemma~\ref{lem:groups_prime_perm}, $\Aut(N)$ is a
subgroup of $\cyc_2^2$. If it is trivial or $\cyc_2$, observe that it can be constructed using (b)
and (c). The only remaining case is when $\Aut(N) \cong \cyc_2^2$. The action of $\cyc_2^2$ on
$V(N)$ can have orbits of sizes $4$, $2$, and $1$. By Orbit-Stabilizer Theorem, each orbit of size
$2$ has also a stabilizer of size $2$, having exactly one non-trivial element. Therefore, there are
at most three types of orbits of size $2$, according to which element of $\cyc_2^2$ stabilizes them.
We give a geometric argument that one of these elements cannot be a stabilizer of an orbit of size
$2$, so there are at most two types of orbits of size $2$.

As argued above, the non-identity elements of $\cyc_2^2$ correspond geometrically to the reflections
$\varphi_v$ and $\varphi_h$ and to the rotation $\rho$; see
Fig.~\ref{fig:action_representations}. The reflection $\varphi_v$ stabilizes those segments which
are parallel to the horizontal axis. The rotation $\rho$ stabilizes those segments which cross the
central point. For both automorphisms, there might be multiple segments stabilized. On the other
hand, the reflection $\varphi_h$ stabilizes at most one segment which lies on the axis of
$\varphi_h$. Further, this segment is stabilized by all elements of $\cyc_2^2$, so it belongs to the
orbit of size $1$. Therefore, there exists no orbit of size $2$ which is stabilized by $\varphi_h$.

Let $G_1$ be the direct product of all $\Aut(T_j)$, one for each orbit of size four.  The groups
$G_2$ and $G_3$ are defined similarly for the orbits of size two stabilized by $\varphi_v$ and
$\rho$, respectively, and $G_4$ for the orbit of size one (if it exists). We have
$$\Aut(T) \cong \bigl (G_1^4 \times G_2^2 \times G_3^2\bigr) \rtimes_{\psi} \cyc_2^2 \times G_4,$$
where $\psi(\varphi_h)$ and $\psi(\varphi_v)$ swap the coordinates as $\varphi_h$ and $\varphi_v$ in
Fig.~\ref{fig:action_representations}. So $\Aut(T)$ can be constructed using (b) and (d).
\end{proof}

\subsection{Automorphism Groups of Bipartite Permutation Graphs.}

We use the modular trees to characterize $\Aut({\rm connected}\ \bipperm)$.  For a connected
bipartite graph, every non-trivial module is an independent set, and the quotient is a prime
bipartite permutation graph.  Therefore, the modular tree $T$ has a prime root node $N$, to which
there are attached leaf nodes which are independent sets.

\begin{figure}[b]
\centering
\includegraphics[scale=0.8]{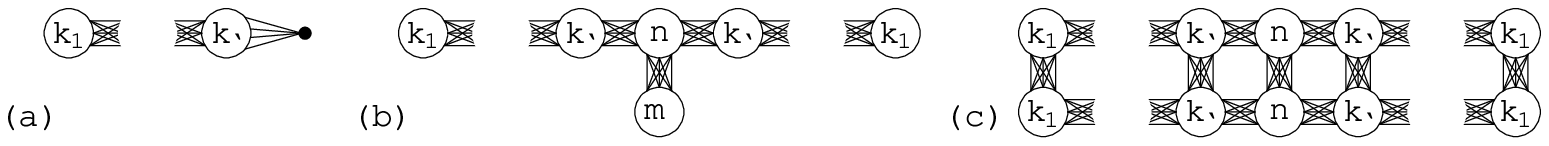}
\caption{Let $G_1 = \sym_{k_1} \times \cdots \times \sym_{k_\ell}$, $G_2 = \sym_n$ and $G_3 =
\sym_m$. The constructed graphs consist of independent sets joined by complete bipartite subgraphs.
They have the following automorphism groups: (a) $G_1$, (b) $G_1 \wr \cyc_2 \times G_2 \times G_3$, (c) $(G_1^4 \times G_2^2) \rtimes
\cyc_2^2$.}
\label{fig:bip_perm_construction}
\end{figure}

\begin{proof}[Proof of Corollary~\ref{cor:aut_groups_bipperm}]
Every abstract group from Corollary~\ref{cor:aut_groups_bipperm} can be constructed as shown in
Fig.~\ref{fig:bip_perm_construction}.  Let $T$ be the modular tree representing $X$.  By
Lemmas~\ref{lem:aut_group_modular_tree} and~\ref{lem:recursive_aut_modular_tree},
$$\Aut(X) \cong \bigl(\Aut(T_1) \times \cdots \times \Aut(T_k)\bigr) \rtimes \Aut(N),$$
where $\Aut(N)$ is isomorphic to a subgroup of $\cyc_2^2$ (by Lemma~\ref{lem:groups_prime_perm}), and
each $\Aut(T_i)$ is a symmetric group since $T_i$ is an independent set.

Consider a permutation representation of $N$ in which the endpoints of the segments, representing
$V(N)$, are placed equidistantly as in Fig.~\ref{fig:action_representations}. By~\cite{bipperm},
there are no segments parallel with the horizontal axis, so the reflections $\varphi_v$ and
$\varphi_h$ fix no segment.  Further, since $N$ is bipartite, there are at most two segments
crossing the central point, so the rotation $\rho$ can fix at most two segments.

\emph{Case 1: $\Aut(N)$ is trivial.} Then $\Aut(X)$ is a direct product of symmetric groups.

\emph{Case 2: $\Aut(N) \cong \cyc_2$.} Let $G_1$ be the direct product of all $\Aut(T_i)$, one for
each orbit of size two. Notice that $\Aut(N)$ is generated by exactly one of $\varphi_v$,
$\varphi_h$, and $\rho$. For $\varphi_v$ or $\varphi_h$, all orbits are of size two, so $\Aut(X)
\cong G_1 \wr \cyc_2$. For $\rho$, there are at most two fixed segments, so $\Aut(X) \cong
G_1 \wr \cyc_2 \times G_2 \times G_3$, where $G_2$ and $G_3$ are isomorphic to $\Aut(T_i)$, for each
of two orbits of size one.

\emph{Case 3: $\Aut(N) \cong \cyc_2^2$.} Then $\Aut(N)$ has no orbits of size $1$, at most one of
size $2$, and all other of size $4$.  Let $G_1$ be the direct product of all $\Aut(T_i)$, one for
each orbit of size $4$, and let $G_2$ be $\Aut(T_i)$ for the orbit of size $2$. We have  $\Aut(X)
\cong (G_1^4 \times G_2^2) \rtimes_{\psi} \cyc_2^2$, where $\psi$ is defined in the proof of
Theorem~\ref{thm:aut_groups_perm}.
\end{proof}

\subsection{$k$-Dimensional Comparability Graphs} \label{sec:kdim}

In this section, we prove that $\Aut(\dim{4})$ contains all abstract finite groups, i.e., each
finite group can be realised as an automorphism group of some $4$-dimensional comparability graph.
Our construction also shows that graph isomorphism testing of $\dim{4}$ is \cGI-complete.  Both
results easily translate to $\dim{k}$ for $k > 4$ since $\dim{4} \subsetneq \dim{k}$.

\heading{The Construction.}
Let $X$ be a graph with $V(X) = \{x_1, \dots, x_n\}$ and $E(X) = \{e_1, \dots,
e_m\}$. We define
$$P = \bigl\{p_i : x_i \in V(X)\bigr\},\qquad
  Q = \{q_{ik} : x_i \in e_k\},\qquad
  R = \bigl\{r_k : e_k \in E(X)\bigr\},$$
where $P$ represents the vertices, $R$ represents the edges and $Q$ represents
the incidences between the vertices and the edges.

The constructed comparability graph $C_X$ is defined as follows, see
Fig.~\ref{fig:example_of_construction}:
$$V(C_X) = P \cup Q \cup R,\qquad E(C_X) = \{p_iq_{ik}, q_{ik}r_k : x_i \in e_k\}.$$
So $C_X$ is created from $X$ by replacing each edge with a path of length four.

\begin{figure}[b]
\centering
\includegraphics[scale=0.8]{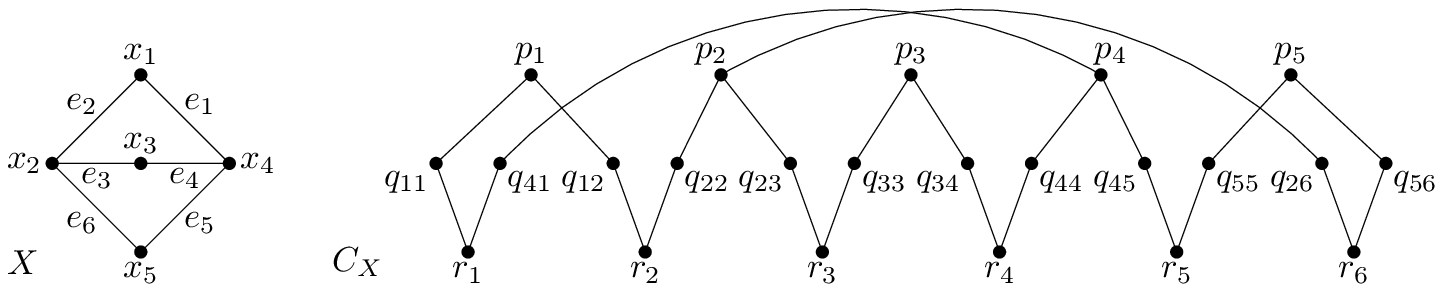}
\caption{The construction $C_X$ for the graph $X = K_{2,3}$.}
\label{fig:example_of_construction}
\end{figure}

\begin{lemma} \label{lem:construction_preserves_aut_group}
For a connected graph $X \not\cong C_n$, $\Aut(C_X) \cong \Aut(X)$.\qed
\end{lemma}


\begin{lemma} \label{lem:bipartite_4dim}
If $X$ is a connected bipartite graph, then $\dm(C_X) \leq 4$.
\end{lemma}

\begin{proof}
We construct four chains such that $L_1 \cap L_2 \cap L_3 \cap L_4$ have two
vertices comparable if and only if they are adjacent in $C_X$.  We describe
linear chains as words containing each vertex of $V(C_X)$ exactly once.  If
$S_1,\dots,S_s$ is a sequence of words, the symbol $\langle S_t :
\mathop\uparrow t \rangle$ is the concatenation $S_1S_2\dots S_s$ and $\langle
S_t : \mathop\downarrow t \rangle$ is the concatenation $S_sS_{s-1}\dots S_1$.
When an arrow is omitted, as in $\langle S_t \rangle$, we concatenate in an
arbitrary order.

First, we define the \emph{incidence string} $I_i$ which codes $p_i$ and its
neighbors $q_{ik}$: $$I_i = p_i \bigl\langle q_{ik} : p_iq_{ik} \in E(C_X)
\bigr\rangle.$$ Notice that the concatenation $I_iI_j$ contains the right edges
but it further contains edges going from $p_i$ and $q_{ik}$ to $p_j$ and
$q_{j\ell}$. We remove these edges by the concatenation $I_jI_i$ in some other
chain. 

Since $X$ is bipartite, let $(A,B)$ be the partition of its vertices. We define
\begin{eqnarray*}
P_A = \{p_i : x_i \in A\},\hskip 0.25em&&\qquad Q_A = \{q_{ik} : x_i \in A\},\\
P_B = \{p_j : x_j \in B\},&&\qquad Q_B = \{q_{jk} : x_j \in B\}.
\end{eqnarray*}
Each vertex $r_k$ has exactly one neighbor in $Q_A$ and exactly one in $Q_B$.

We construct the four chains as follows:
\begin{eqnarray*}
L_1 &=& \langle p_i : p_i \in P_A \rangle \langle r_kq_{ik} : q_{ik} \in Q_A, \mathop\uparrow k\rangle
\langle I_j : p_j \in P_B, \mathop\uparrow j\rangle,\\
L_2 &=& \langle p_i : p_i \in P_A \rangle \langle r_kq_{ik} : q_{ik} \in Q_A, \mathop\downarrow k\rangle
\langle I_j : p_j \in P_B, \mathop\downarrow j\rangle,\\
L_3 &=& \langle p_j : p_j \in P_B \rangle \langle r_kq_{jk} : q_{jk} \in Q_B, \mathop\uparrow k\rangle
\langle I_i : p_i \in P_A, \mathop\uparrow i\rangle,\\
L_4 &=& \langle p_j : p_j \in P_B \rangle \langle r_kq_{jk} : q_{jk} \in Q_B, \mathop\downarrow k\rangle
\langle I_i : p_i \in P_A, \mathop\downarrow i\rangle.
\end{eqnarray*}

\begin{figure}[t]
\centering
\includegraphics[scale=0.8]{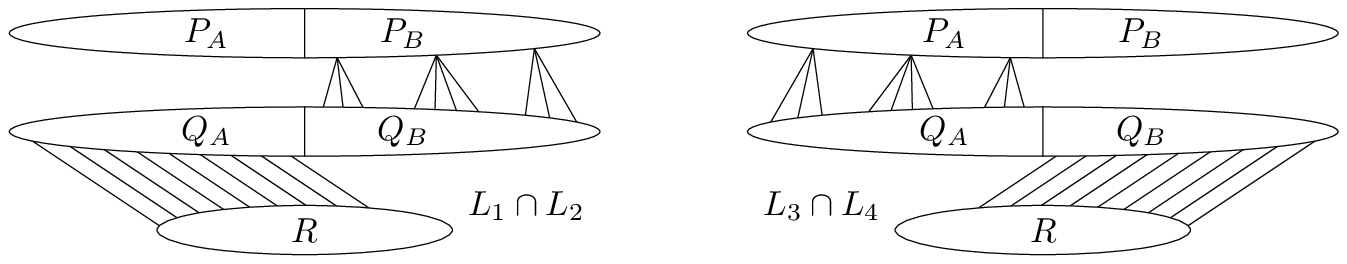}
\caption{The forced edges in $L_1 \cap L_2$ and $L_3 \cap L_4$.}
\label{fig:proof_idea}
\end{figure}

The four defined chains have the following properties, see Fig.~\ref{fig:proof_idea}:
\begin{itemize}
\item
The intersection $L_1 \cap L_2$ forces the correct edges between $Q_A$ and $R$
and between $P_B$ and $Q_B$. It poses no restrictions between $Q_B$ and $R$ and
between $P_A$ and the rest of the graph.
\item
Similarly the intersection $L_3 \cap L_4$ forces the correct edges between $Q_B$
and $R$ and between $P_A$ and $Q_A$. It poses no restrictions between $Q_A$ and
$R$ and between $P_B$ and the rest of the graph.
\end{itemize}
It is routine to verify that the intersection $L_1 \cap L_2 \cap L_3 \cap L_4$ is correct.

\medskip
\noindent\emph{Claim 1: The edges in $Q \cup R$ are correct.} For every $k$, we
get $r_k$ adjacent to both $q_{ik}$ and $q_{jk}$ since it appear on the left in
$L_1,\dots,L_4$. On the other hand, $q_{ik}q_{jk} \notin E(C_X)$ since they are
ordered differently in $L_1$ and $L_3$.

For every $k < \ell$, there are no edges between
$N[r_k] = \{r_k,q_{ik},q_{jk}\}$ and $N[r_\ell] = \{r_\ell,q_{s\ell},q_{t\ell}\}$. This can be shown
by checking the four orderings of these six elements:
\begin{eqnarray*}
\text{in $L_1$:}\quad r_k q_{ik}\, \boxed{r_\ell q_{s\ell}}\, q_{jk}\, \boxed{q_{t \ell}},\quad&&\quad
\text{in $L_2$:}\quad \boxed{r_\ell q_{s\ell}}\, r_k q_{ik} q_{jk}\, \boxed{q_{t \ell}},\\
\text{in $L_3$:}\quad r_k q_{jk}\, \boxed{r_\ell q_{t\ell}}\, q_{ik}\, \boxed{q_{s \ell}},\quad&&\quad
\text{in $L_4$:}\quad \boxed{r_\ell q_{t\ell}}\, r_k q_{jk} q_{ik}\, \boxed{q_{s \ell}},
\end{eqnarray*}
where the elements of $N[r_\ell]$ are boxed.\claimqed

\medskip
\noindent\emph{Claim 2: The edges in $P$ are correct.} We show that there are no
edges between $p_i$ and $p_j$ for $i \ne j$ as follows.  If both belong to $P_A$
(respectively, $P_B$), then they are ordered differently in $L_3$ and $L_4$
(respectively, $L_1$ and $L_2$). If one belongs to $P_A$ and the other one to
$P_B$, then they are ordered differently in $L_1$ and $L_3$.\claimqed

\medskip
\noindent\emph{Claim 3: The edges between $P$ and $Q \cup R$ are correct.} For
every $p_i \in P$ and $r_k \in R$, we have $p_ir_k \notin E(C_X)$ because they
are ordered differently in $L_1$ and $L_3$. On the other hand, $p_iq_{ik} \in
E(C_X)$, because $p_i$ is before $q_{ik}$ in $I_i$, and for $p_i \in P_A$ in
$L_1$ and $L_2$, and for $p_i \in P_B$ in $L_3$ and $L_4$.

It remains to show that $p_iq_{jk} \notin E(C_X)$ for $i \ne j$.  If both $p_i$
and $p_j$ belong to $P_A$ (respectively, $P_B$), then $p_i$ and $q_{jk}$ are
ordered differently in $L_3$ and $L_4$ (respectively, $L_1$ and $L_2$). And if
one belongs to $P_A$ and the other one to $P_B$, then $p_i$ and $q_{jk}$ are
ordered differently in $L_1$ and $L_3$.\claimqed
\medskip

These three established claims show that comparable pairs in the intersection $L_1 \cap L_2
\cap L_3 \cap L_4$ are exactly the edges of $C_X$, so $C_X \in \dim4$.
\end{proof}

\heading{Universality of \bolddim{k}.} We are ready to prove Theorem~\ref{thm:kdim_aut_groups_and_gi}.

\begin{proof}[Proof of Theorem~\ref{thm:kdim_aut_groups_and_gi}]
We prove the statement for \dim4. Let $X$ be a connected graph such that $X
\not\cong C_n$.  First, we construct the bipartite incidence graph $Y$ between $V(X)$ and $E(X)$.
Next, we construct $C_Y$ from $Y$. From
Lemma~\ref{lem:construction_preserves_aut_group} it follows that $\Aut(C_Y) \cong \Aut(Y) \cong
\Aut(X)$ and by Lemma~\ref{lem:bipartite_4dim}, we have that $C_Y \in \dim 4$.  Similarly, if two
graphs $X_1$ and $X_2$ are given, we construct $C_{Y_1}$ and $C_{Y_2}$ such that $X_1 \cong X_2$ if
and only if $C_{Y_1} \cong C_{Y_2}$; this polynomial-time reduction shows \cGI-completeness of graph
isomorphism testing.

The constructed graph $C_Y$ is a prime graph. We fix the transitive orientation in which $P$ and $R$
are the minimal elements and get the poset $P_Y$ with $\Aut(P_Y) \cong \Aut(C_Y)$.  Hence, our
results translate to posets of the dimension at most four.
\end{proof}

\section{Algorithms for Computing Automorphism Groups} \label{sec:algorithms}

Using PQ-trees, Colbourn and Booth~\cite{autmorphism_algorithms_trees_int_planar} give a linear-time
algorithm to compute permutation generators of the automorphism group of an interval graph. We are
not aware of any such algorithm for circle and permutation graphs, but some of our results might be
known from the study of graph isomorphism problem~\cite{hsu1995m,permutation_isomorphism}. 

We briefly explain algorithmic implications of our results which allow to compute automorphism
groups of studied classes in terms of $\cyc_n$, $\dih_n$ and $\sym_n$, and their group products.
This description is better than just a list of permutations generating $\Aut(X)$. Many tools of the
computational group theory are devoted to getting better understanding of an unknown group.  Our
description gives this structural understanding of $\Aut(X)$ for free.

For interval graphs, we get a linear-time algorithm by computing an
MPQ-tree~\cite{incremental_linear_int_recognition} and finding its symmetries.  For circle graphs,
our description easily leads to a polynomial-time algorithm, by computing the split tree for each
connected component and understanding its symmetries. The best algorithm for computing split trees
runs in almost linear time~\cite{gioan2013practical}. For permutation graph, we get a linear-time
algorithm by computing the modular decomposition~\cite{mcconnell_spinrad} and finding symmetries of
prime permutation graphs.

\section{Open Problems} \label{sec:conclusions}

We conclude this paper with several open problems concerning automorphism groups of other
intersection-defined classes of graphs; for an overview see~\cite{agt,egr}.

\emph{Circular-arc graphs} (\ca) are intersection graphs of circular arcs and they naturally
generalize interval graphs. Surprisingly, this class is very complex and quite different from
interval graphs. Hsu~\cite{hsu1995m} relates them to circle graphs.

\begin{problem}
What is $\Aut(\ca)$?
\end{problem}

Let $Y$ be any fixed graph. The class $\graphs{Y}$ consists of all intersections graphs of connected
subgraphs of a subdivision of $Y$. Observe that $\graphs{K_2} = \int$ and we have an infinite
hierarchy between $\int$ and $\chor$ is formed by $\graphs{T}$ for a tree $T$, for which $\int
\subseteq \graphs{T} \subsetneq \chor$. If $Y$ contains a cycle, then $\graphs{Y} \not\subseteq
\chor$. For instance, $\graphs{K_3} = \ca$.

\begin{conjecture}
For every fixed graph $Y$, the class $\graphs{Y}$ is non-universal.
\end{conjecture}

The last open problem involves the open case of \dim3.

\begin{conjecture}
The class \dim3 is universal and its graph isomorphism problem is \cGI-complete.
\end{conjecture}

\heading{Acknowlegment.} We would like to thank to Roman Nedela for many comments. A part of this
work was done during a visit at Matej Bel University. 

\bibliographystyle{amsplain}
\bibliography{aut_groups_of_geom_rep_graphs}

\end{document}



\section{}
\subsection{}

\begin{theorem}[Optional addition to theorem head]
\end{theorem}

\begin{proof}[Optional replacement proof heading]
\end{proof}

\begin{figure}
\includegraphics{filename}
\caption{text of caption}
\label{}
\end{figure}


\begin{equation}
\end{equation}

\begin{equation*}
\end{equation*}

\begin{align}
  &  \\
  &
\end{align}
